\documentclass{msripub}
\let\newcommand\providecommand
\usepackage{amsmath}
\usepackage{book51}
\usepackage{msribib}
\usepackage{mathrsfs}
\usepackage{xypic}
\let\mathcal\mathscr

\cleardoublepage

\let\iff\Leftrightarrow
\let\obf\subsubsectionfont
\def\div{\mskip1.5mu|\mskip1.5mu}
\let\implies\Rightarrow

\newcommand{\also}{\quad\text{and}\quad}
\newcommand{\qnd}{\quad\text{and}\quad} 

\newcommand{\foral}{\quad\text{for all}\ }

\newcommand{\where}{\quad\text{where}\ }

\newcommand{\even}{{\operatorname{even}}}
\newcommand{\odd}{{\operatorname{odd}}}

\newcommand{\sgn}[1]{(-1)^{|#1|}}
\newcommand{\sign}[2]{(-1)^{|#1||#2|}}

\newcommand{\col}{\colon}
\newcommand{\dd}{\partial}

\newcommand{\nat}{{}^{\natural}}
\newcommand{\wt}{\widetilde}

\newcommand{\alto}{\xleftarrow}
\newcommand{\arto}{\xrightarrow}

\newcommand{\tra}{\twoheadrightarrow}

\newcommand{\id}{\operatorname{id}}

\newcommand{\Ker}{\operatorname{Ker}}

\newcommand{\aug}[1]{\operatorname{I}(#1)}

\newcommand{\der}[2]{\operatorname{Der}(#1\hskip.8pt;#2)}
\newcommand{\inn}[2]{\operatorname{IDer}(#1\hskip.8pt;#2)}
\newcommand{\dual}[1]{{#1}{}^*}

\newcommand{\gring}[2]{{#1}[#2]}

\newcommand{\rank}{\operatorname{rank}}

\newcommand{\sym}{\operatorname{Sym}}
\newcommand{\soc}[1]{\operatorname{socle}(#1)}

\newcommand{\fdim}{\operatorname{fdim}} 
\newcommand{\pdim}{\operatorname{pdim}}
\newcommand{\injdim}{\operatorname{injdim}}

\newcommand{\gc}[3]{\operatorname{H}^{#1}(#2\hskip.8pt,#3)}

\newcommand{\EC}[3]{\operatorname{E}_{#1}^{#2,#3}}

\newcommand{\Ext}[4]{{\operatorname{Ext}^{#1}_{#2}(#3,#4)}}

\newcommand{\Hom}[3]{\operatorname{Hom}_{#1}(#2,#3){}}
\newcommand{\tensor}[3]{{#1}\otimes_{#2}{#3}}
\newcommand{\Tor}[4]{\operatorname{Tor}_{#1}^{#2}(#3,#4){}}

\newcommand{\BZ}{{\mathbb Z}}

\newcommand{\bsP}{{\mathbf P}}

\newcommand{\fm}{{\mathfrak m}}
\newcommand{\fn}{{\mathfrak n}}

\newcommand{\eps}{{\epsilon}}

\newcommand{\vf}{{\varphi}}

\newtheorem[{}{\it\specialdigits}]{theorem}{Theorem}[section]
\newtheorem[{}{\it\specialdigits}]{Theorem}{Theorem}\unnumbered{Theorem}
\newtheorem[{}{\it\specialdigits}]{maschke}[Theorem]{Maschke's theorem}
\newtheorem[{}{\it\specialdigits}]{proposition}[theorem]{Proposition}
\newtheorem[{}{\it\specialdigits}]{Proposition}[Theorem]{Proposition}
\newtheorem[{}{\it\specialdigits}]{lemma}[theorem]{Lemma}
\newtheorem[{}{\it\specialdigits}]{Lemma}[Theorem]{Lemma}
\newtheorem[{}{\it\specialdigits}]{corollary}[theorem]{Corollary}
\newtheorem[{}{\it\specialdigits}]{Corollary}[Theorem]{Corollary}
\newtheorem[{}{\it\specialdigits}]{property}[theorem]{}
\newtheorem[{}{\it\specialdigits}]{Property}[Theorem]{Property}
\newtheorem[{}{\it\specialdigits}]{noname}[theorem]{}

\newtheorem{exercise}[theorem]{Exercise}
\newtheorem{Exercise}[Theorem]{Exercise}

\newtheorem{chunk}[theorem]{}
\makeatletter
\def\thrmintro#1#2{\expandafter\let\expandafter\@tempa#2\ifx\@tempa\relax
\else(#2) \fi\def\@tempa{#1}\ifx\@tempa\@empty\expandafter\@gobble\else#1\fi}
\makeatother

\newtheorem{ramble}[Theorem]{Ramble}

\newtheorem{remark}[theorem]{Remark}

\newtheorem{case}[Theorem]{Case}

\title
[Modules and cohomology over group algebras]
{Modules and Cohomology over Group Algebras:
One Commutative Algebraist's Perspective}

\author{Srikanth Iyengar}

\address{305  Avery Hall\\
Department of Mathematics\\
University of Nebraska\\
Lincoln, NE 68588\\
United States}

\email{iyengar@math.unl.edu}

\thanks{Part of this article was written while the author was funded by 
a grant from the NSF}


\begin{document}

\begin{abstract}
This article explains basic constructions and results on group
algebras and their cohomology, starting from the point of view of
commutative algebra. It provides the background necessary for a novice
in this subject to begin reading Dave Benson's article in this volume. 
\end{abstract}

\maketitle

\tableofcontents

\section*{Introduction}

The available
accounts of group algebras and group cohomology [\citeNP{Benson:1991a};
  \citeyearNP{Benson:1991b}; \citeNP{Brown:1982a,Evens:1991a}]
 are all written for the mathematician on the
street.
This one is written for commutative algebraists by one of their own. There is a point to
such an exercise: though group algebras are typically noncommutative, module theory over
them shares many properties with that over commutative rings. Thus, an exposition that
draws on these parallels could benefit an algebraist familiar with the commutative world.
However, such an endeavour is not without its pitfalls, for often there are subtle
differences between the two situations.  I have tried to draw attention to similarities
and to discrepancies between the two subjects in a series of commentaries on the text that
appear under the rubric Ramble\footnote{ This word has at least two meanings: ``a
  leisurely walk'', or ``to talk or write in a discursive, aimless way''; you can decide
  which applies. By the by, its etymology, at least according to
  www.dictionary.com, might amuse you.}.

The approach I have adopted toward group cohomology is entirely algebraic. However, one
cannot go too far into it without some familiarity with algebraic topology. To gain an
appreciation of the connections between these two subjects, and for a history of group
cohomology, one might read \cite{Be:hist,Mc:hist}.

In preparing this article, I had the good fortune of having innumerable `chalk-and-board'
conversations with Lucho Avramov and Dave Benson.  My thanks to them for all these, and to
the Mathematical Sciences Research Institute for giving me an opportunity to share a roof
with them, and many others, during the Spring of 2003. It is also a pleasure to thank
Kasper Andersen, Graham Leuschke, and Claudia Miller for their remarks and suggestions.

\section{The Group Algebra}
\label{The group algebra}

Let $G$ be a group, with identity element $1$, and let $k$ be a field.  Much of what is
said in this section is valid, with suitable modifications, more generally when $k$ is a
commutative ring.  Let $\gring kG$ denote the $k$-vector space with basis the elements of
$G$; thus $\gring kG = \bigoplus_{g\in G}kg$.  The product on $G$ extends to an
associative multiplication on $\gring kG$: for basis elements $g$ and $h$, one has $g\cdot
h=gh$, where the product on the right is taken in $G$, while the product of arbitrary
elements is specified by the distributive law and the rule $a\cdot g = g\cdot a$ for
$a\in k$.  The identity element $1$ is the identity in $\gring kG$.  The $k$-linear ring
homomorphism $\eta\col k\to\gring kG$ with $\eta(1)=1$ makes $\gring kG$ a $k$-algebra.
This is the \emph{group algebra} of $G$ with coefficients in $k$.

Note that $\gring kG$ is commutative if and only if the group $G$ is abelian.  Moreover,
it is finite-dimensional as a $k$-vector space precisely when $G$ is finite.
  
An important part of the structure on $\gring kG$ is the augmentation of $k$-algebras
$\eps\col \gring kG\to k$ defined by $\eps(g)=1$ for each $g\in G$. Through $\eps$ one can
view $k$ as a $\gring kG$-bimodule.  The kernel of $\eps$, denoted $\aug G$, is the
$k$-subspace of $\gring kG$ with basis $\thickmuskip 5mu minus 3mu \{g{-}1\mid g\in G\}$; it is a two-sided ideal,
called the \textit{augmentation ideal} of~$G$.  For every pair of elements $g,h$ in $G$, 
the following relations hold in the group algebra:
\begin{gather*}
g^{-1}-1 = g^{-1}(1-g),\\
gh - 1 = g(h-1) + (g-1) = (g-1)h + (h-1).
\end{gather*}
Thus, if a subset $\{g_\lambda\}_{\lambda\in\Lambda}$ of $G$, with $\Lambda$ an index set,
generates the group, the subset $\{g_{\lambda}-1\}_{\lambda\in\Lambda}$ of $\gring
kG$ generates $\aug G$ both as a left ideal and as a right ideal.{\looseness=-1\par}

\begin{chunk}{\obf Functoriality.}
\label{gring:functor}
The construction of the group algebra is functorial: given a group homomorphism $\vf\col
G_1\to G_2$, the $k$-linear map
\[
\gring k{\vf}\col \gring k{G_1}\to \gring k{G_2}, \where g\mapsto \vf(g),
\]
is a homomorphism of $k$-algebras, compatible with augmentations. Its kernel is generated
both as a left ideal and as a right ideal by the set $\{g-1\mid g\in \Ker\vf\}$.

For example, when $N$ is a normal subgroup of a group $G$, the canonical surjection
$G\to G/N$ induces the surjection of $k$-algebras $\gring kG\to\gring k{G/N}$. Since its
kernel is generated by the set $\{n-1\mid n\in N\}$, there is a natural isomorphism of
$k$-algebras
\[
\gring k{G/N}\cong k\otimes_{\gring kN}\gring kG = \frac{\gring kG}{\aug N\gring kG}.
\]
\end{chunk}

Let me illustrate these ideas on a few simple examples.

\begin{chunk}{\obf Cyclic groups.}
\label{gring:cyclic}
The group algebra of the infinite cyclic group is $k[x^{\pm 1}]$, the algebra of Laurent
polynomials in the variable $x$. Here $x$ is a generator of the group; its inverse is
$x^{-1}$. The augmentation maps $x$ to $1$, and the augmentation ideal is generated, as an
ideal, by $x-1$.

In view of \eqref{gring:functor}, the group algebra of the cyclic group of order $d$ is
$k[x]/(x^d-1)$, and the augmentation ideal is again generated by $x-1$.
\end{chunk}

\begin{chunk}{\obf Products of groups.}
\label{gring:products}
Let $G_1$ and $G_2$ be groups. By \eqref{gring:functor}, for $n=1,2$ the canonical
inclusions $\iota_n\col G_n\to G_1\times G_2$ induce homomorphisms of $k$-algebras $\gring
k{\iota_n}\col \gring k{G_n}\to \gring k{G_1\times G_2}$. Since the elements in the image
of $\gring k{\iota_1}$ commute with those in the image of $\gring k{\iota_2}$, one obtains
a homomorphism of augmented $k$-algebras
\begin{align*}
{\gring k{G_1}}&\otimes_k{\gring k{G_2}} \to \gring k{G_1\times G_2},\\
            g_1&\otimes_k\rlap{$g_2$}\phantom{\gring k{G_2}} \mapsto (g_1,g_2).
\end{align*}
This is an isomorphism since it maps the basis $\{g_1\otimes_kg_2\mid g_i\in G_i\}$ of
the $k$-vector space $\gring k{G_1}\otimes_k\gring k{G_2}$ bijectively to the basis
$\{(g_1,g_2)\mid g_i\in G_i\}$ of $\gring k{G_1\times G_2}$.  For this reason, the group
algebra of $G_1\times G_2$ is usually identified with $\tensor{\gring k{G_1}}k{\gring
  k{G_2}}$.
\end{chunk}

\begin{chunk}{\obf Abelian groups.}
\label{gring:abelian}
Let $G$ be a finitely generated abelian group. The structure theorem for such groups tells
us that there are nonnegative numbers $n$ and $d_1,\dots,d_m$, with $d_j\geq 2$ and
$d_{i+1}\div d_i$, such that
\[ 
G = \BZ^n \oplus\frac \BZ{(d_1\BZ)}\oplus\cdots \oplus\frac \BZ{(d_m\BZ)}.
\]
The description of the group algebra of cyclic groups given in \eqref{gring:cyclic}, in
conjunction with the discussion in \eqref{gring:products}, yields
\[
\gring kG= \frac{k[x_1^{\pm 1},\dots,x_n^{\pm
    1},y_1,\dots,y_m]}{(y_1^{d_1}{-}1,\dots,y_m^{d_m}{-}1)}
\]
The augmentation is given by $x_i\mapsto 1$ and $y_j\mapsto 1$, the augmentation ideal is
generated by $\{x_1{-}1,\dots,x_n{-}1,y_1{-}1,\dots,y_m{-}1\}$.
\end{chunk}

\begin{ramble}
Observe: the group algebra in \eqref{gring:abelian} above is a complete intersection.
\end{ramble}

\begin{chunk}{\obf Finite $p$-groups.}
\label{gring:pgroups}
Let $R$ be a ring; it need not be commutative. Recall that the intersection of all its
left maximal ideals is equal to the intersection of all its right maximal ideals, and
called the Jacobson radical of $R$. Thus, $R$ has a unique left maximal ideal exactly when
it has a unique right maximal ideal, and then these ideals coincide.  In this case, one
says that $R$ is \textit{local}; note that the corresponding residue ring is a division
ring; for details see \cite[XVII~\S6]{La}, for example.
  
Suppose that the characteristic of $k$ is $p$, with $p\geq 2$. Let $G$ be a finite
\textit{$p$-group}, so that the order of $G$ is a power of $p$.  I claim:

\medskip

\textit{The  group algebra $\gring kG$ is local with maximal ideal $\aug G$.} 

\medskip

Indeed, it suffices to prove (and the claim is equivalent to): the augmentation ideal
$\aug G$ is nilpotent. Now, since $G$ is a $p$-group, its centre $Z$ is nontrivial, so
\eqref{gring:functor} yields an isomorphism of $k$-algebras
\[
\frac{\gring kG}{\aug Z \gring kG} \cong \gring k{G/Z}.
\]
Since the order of $G/Z$ is strictly less than that of $G$, one can assume that $\aug
{G/Z}$ is nilpotent.  By the isomorphism above, this entails ${\aug G}^n\subseteq \aug Z
\gring kG$, for some positive integer $n$. Now $Z$ is an abelian $p$-group, so $\aug Z$ is
nilpotent, by \eqref{gring:abelian}. Since $\aug Z$ is in the centre of $\gring kG$,
one obtains that $\aug G$ is nilpotent, as claimed.
 \end{chunk}

The converse also holds:

\begin{exercise}
\label{local=pgroup}
Let $G$ be a finite group and $p$ the characteristic of $k$.  Prove that if the ring
$\gring kG$ is local, then $G$ is a $p$-group. (Hint: $\gring kG$ has finite rank
over $k$, so its nilradical is equal to its Jacobson radical.)
 \end{exercise}

\begin{chunk}{\obf The diagonal map.}
\label{gring:diagonal}
Let $G$ be a group and let $G\to G\times G$ be the diagonal homomorphism, given by $g\mapsto
(g,g)$. Following \eqref{gring:products}, one identifies the group ring of $G\times G$
with $\tensor{\gring kG}k{\gring kG}$, and then the diagonal homomorphism induces a
homomorphism of augmented $k$-algebras
\[
\Delta \col \gring kG \to \tensor{\gring kG}k{\gring kG}, \where g\mapsto \tensor gkg.
\]
This is called the \textit{diagonal} homomorphism, or \textit{coproduct}, of the group algebra
$\gring kG$.

There is another piece of structure on the group algebra: the map $G\to G$ given by
$g\mapsto g^{-1}$ is an anti-isomorphism of groups, and hence induces an anti-isomorphism
of group algebras 
\[
\sigma\col \gring kG\to \gring kG,
\]
that is to say, $\sigma$ is an isomorphism of additive groups with
$\sigma(rs)=\sigma(s)\sigma(r)$. The map $\sigma$ is referred to as the \textit{antipode} of
the group algebra. It commutes with the diagonal map, in the sense that
\[
\sigma^{(G\times G)} \circ \Delta^G = \Delta^G\circ \sigma^G.
\]
Here are the salient properties of the diagonal and the antipode:
\begin{enumerate}
\item[{\rm(a)}] $\Delta$ is a homomorphism of augmented $k$-algebras;
\item[{\rm(b)}] $\Delta$ is co-associative, in that the following diagram commutes:
\[
\xymatrixrowsep{2pc} \xymatrixcolsep{3pc} \xymatrix{ \gring kG \ar@{->}[r]^{\Delta}
  \ar[d]^{\Delta} &
  \gring kG\otimes_k\gring kG\ar[d]^{\Delta\otimes_k1} \\
  \gring kG\otimes_k\gring kG\ar[r]^-{1\otimes_k\Delta} &\gring kG\otimes_k\gring
  kG\otimes_k\gring kG }
\]
\item[{\rm(c)}] The following diagram commutes:
\[
\xymatrixrowsep{2pc} \xymatrixcolsep{3pc} \xymatrix{
  &\gring kG\ar[dl]_{\cong}\ar[d]^{\Delta}\ar[dr]^{\cong}& \\
  k\otimes_k \gring kG\ar@{<-}[r]^{\eps\otimes_k 1} &\gring kG\otimes_k\gring kG &
  \ar@{<-}[l]_{1\otimes_k\eps}\gring kG\otimes_kk}
 \]
 This property is paraphrased as: $\eps$ is a co-unit for $\Delta$.
\item[{\rm(d)}] For each element $r\in \gring kG$, if $\Delta(r)=\sum_{i=1}^n (r_i'\otimes_kr_i'')$,
  then
\[
\sum_{i=1}^n\sigma(r_i')r_i'' = \eta(\eps(r))= \sum_{i=1}^n r_i'\sigma(r_i'')
\]
\end{enumerate}
\end{chunk}

Taking these properties as the starting point, one arrives at the following notion.

\begin{chunk}{\obf Hopf algebras.}
\label{hopf algebras}
An augmented $k$-algebra $H$, with unit $\eta\col k\to H$ and augmentation $\eps\col H\to
k$ with $k$-linear homomorphisms $\Delta\col H\to H\otimes_kH$ and $\sigma\col H\to H$
satisfying conditions (a)--(d) listed above, is said to be a \textit{Hopf algebra}.  Among
these, (b) and (c) are the defining properties of a \textit{coalgebra} with diagonal
$\Delta$; see \cite{Mo} or \cite{Sw}.  Property (a) says that the algebra and coalgebra
structures are compatible.  At first\emdash and perhaps second and third\emdash glance, property (d)
appears mysterious. Here is one explanation that appeals to me: The diagonal
homomorphism endows the $k$-vector space $\Hom kHH$ with the structure of a $k$-algebra,
with the product of elements $f$ and $g$  given by
\[
(f\star g)(r) = \sum_{i=1}^n f(r_i')g(r_i''), \where \Delta(r)=\sum_{i=1}^n
(r_i'\otimes_kr_i'').
\]
This is called the \textit{convolution product} on $\Hom kHH$; its unit is the element
$\eta\circ\eps$. Condition (d) asserts that $\sigma$ is the inverse of
the identity on $H$.

The group algebra is the prototypical example of a Hopf algebra, and many constructions
and results pertaining to them are best viewed in that generality;
see \cite[Chapter~3]{Benson:1991a}.  There is another good source of Hopf algebras, close to
home: the coordinate rings of algebraic groups. You might, as I did, find it entertaining
and illuminating to write down the Hopf structure on the coordinate ring of the circle
$x^2+y^2=1$.

If this all too brief foray into Hopf algebras has piqued your curiosity and you wish to
know more, you could start by reading Bergman's charming introduction 
\cite{Be}; if you prefer to jump right into the thick of things, then \cite{Mo} is the
one for you.
 \end{chunk}

\section{Modules over Group Algebras}
\label{sec:modules}

This section is an introduction to modules over group algebras.  
When $G$ is a \textit{finite} group, the $k$-algebra $\gring kG$ is
finite-dimensional, that is to say, of finite rank over $k$. Much of
the basic theory for modules 
over finite group algebras is just a specialization of the theory for finite-dimensional
algebras. For example, I hinted in Exercise \eqref{local=pgroup} that for finite group
algebras, the nilradical coincides with the Jacobson radical; this holds, more generally,
for any finite-dimensional $k$-algebra. Here I will focus on two crucial
concepts: the Jordan--H\"older theorem and the Krull--Schmidt property.

\begin{chunk}{\obf The Jordan--H\"older theorem.}
\label{jordan-holder}
Let $R$ be a ring and $M$ an $R$-module.  It is clear that when $M$ is both artinian and
noetherian it has a \textit{composition series}: a series of submodules $0= M_l\subset
M_{l-1}\subset \cdots \subset M_1\subset M_0=M$ with the property that the subfactors
$M_i/M_{i+1}$ are \textit{simple}, that is to say, they have no proper submodules. It turns
out that if $0= M'_{l'}\subset M'_{l'-1}\subset \cdots \subset M'_1\subset M'_0=M$ is
another composition series, then $l=l'$ and, for $1\leq i,j\leq l$, the factors
$M_i/M_{i-1}$ are a permutation of the factors $M'_j/M'_{j-1}$.  This is a consequence of
the Jordan--H\"older theorem, which says that for each $R$-module, any two series (not
necessarily composition series) of submodules can be refined to series of the same length
and with the same subfactors.
  
Suppose that $R$ is artinian; for example, $R$ may be a finite-dimensional $k$-algebra,
or, more specifically, a finite group algebra. In this case every finite, by which I mean
`finitely generated', module over it is both artinian and noetherian and so has a
composition series. Here is one consequence: since every simple module is a quotient of
$R$, all the simple modules appear in a composition series for $R$, and so there can only
be finitely many of them.
\end{chunk}

\begin{chunk}{\obf Indecomposable modules.}
\label{indecomposables}
Recall that a module is said to be \textit{indecomposable} if it has no nontrivial direct
summands. It is clear that a simple module is indecomposable, but an indecomposable module
may be far from simple\emdash in either sense of the word.  For example, over a commutative ring, the
only simple modules are the residue fields, whereas it is usually not possible to classify
all the indecomposable modules; I will pick up on this point a few paragraphs down the
road.  For now, here are a couple of remarks that are useful to keep in
mind when dealing with indecomposability; see the discussion in
\eqref{modules:klein}.
 \end{chunk}
 
In this sequel, when I say $(R,\fm,k)$ is a local ring, I mean that $R$ is local, with
maximal ideal $\fm$ and residue ring $k$.

\begin{exercise}
\label{indecomposable:socleM}
Let $(R,\fm,k)$ be a commutative local ring. Prove that if $M$ is indecomposable, then $\soc M
\subseteq \fm M$.
\end{exercise}

\begin{exercise}
\label{indecomposable:socleR}
Let $R$ be a commutative local Gorenstein ring and $M$ an indecomposable $R$-module. Prove
that if $\soc R\cdot M\ne 0$, then $M\cong R$.
\end{exercise}

\begin{chunk}{\obf The Krull--Schmidt property.}
\label{krull-schimdt}
Let $R$ be a ring. It is not hard to see that each finite $R$-module can be broken down
into a finite direct sum of indecomposables.  The ring $R$ has the \textit{Krull--Schmidt
  property} if for each finite $R$-module such a decomposition is unique up to a
permutation of the indecomposable factors: if
 \[
\bigoplus_{i=1}^m M_i  \ \cong \  \bigoplus_{j=1}^n N_j,
 \]
with each $M_i$ and $N_j$ indecomposable, then $m=n$, and, with a possible rearrangement
of the $N_j$, one has $M_i\cong N_i$ for each $i$.

For example, complete commutative noetherian local rings have this property; see
\cite[(2.22)]{Swa}. In the present context, the relevant result is that artinian
rings have the Krull--Schmidt property \cite[(1.4.6)]{Benson:1991a}.  When $G$ is a finite group,
$\gring kG$ is artinian; in particular, it has the Krull--Schmidt property.
\end{chunk}

The Krull--Schmidt property is of great help in studying modules over group algebras, for it
allows one to focus on the indecomposables. The natural question arises: when does the
group algebra have only finitely many isomorphism classes of indecomposable modules? In
other words, when is the group algebra of \textit{finite representation type}?  This is the
case, for example, when every indecomposable module is simple, for there are only finitely
many of them; see \eqref{jordan-holder}. There is an important context when this happens:
when the characteristic of $k$ is coprime to the order of the group. This is a consequence
of Maschke's Theorem:

\begin{theorem}[Maschke]
\label{maschke:real}
Let $G$ be a finite group such that $|G|$ is coprime to the the characteristic of $k$.  Each
short exact sequence of $\gring kG$-modules splits.
\end{theorem}
\begin{proof}
  Let $0\to L\to M\arto{\pi} N\to 0$ be an exact sequence of $\gring kG$-modules. Since
  $k$ is a field, $\pi$ admits a $k$-linear section; let $\sigma\col N \to M$ be one such.
  It is not hard to verify that the map
\[
\wt\sigma \col N \to M,\quad \where \wt\sigma(n) = \frac 1{|G|}\,
\sum_{g\in G} g\sigma(g^{-1}n) \foral n\in N,
\]
is $\gring kG$-linear, and that $\pi\circ \wt\sigma =\id^N$. Thus, the exact sequence
splits, as desired.
\end{proof}

This theorem has a perfect converse: if each short exact sequence of $\gring kG$-modules
splits, the characteristic of $k$ is coprime to $|G|$. In fact, it suffices that the
exact sequence $0\to \aug G \to \gring kG\arto{\eps} k\to 0$ splits. The proof is
elementary, and is recommended as an exercise; I will offer a solution in the proof of
Theorem \eqref{maschke}.{\looseness=-1\par}

A group algebra can have finite representation type even if not every indecomposable
module is simple:

\begin{chunk}{\obf Finite cyclic groups.}  
\label{modules:cyclic}
In describing this example, it is convenient to let $p$ denote $1$ when the
characteristic of $k$ is $0$, and the characteristic of $k$ otherwise.  

Let $G$ be a finite cyclic group. Write $|G|$ as $p^nq$, where $n$ is a nonnegative
integer and $p$ and $q$ are coprime. Let $R=k[x]/(x^{{p^n}q}-1)$, the group algebra.  The
binomial theorem in characteristic $p$ yields $x^{p^{n}q}-1 = (x^q-1)^{p^n}$,
so the Jacobson radical of $R$ is $(x^q-1)$.  In $k[x]$, the polynomial $x^q-1$ breaks up
into a product of distinct irreducible polynomials:
\[
x^q - 1 = \prod_{i=1}^df_i(x),\quad \text{with} \ \ \sum_{i=1}^d\deg(f_i(x))=q.
\]
Since the ideals $(f_i(x)^{p^n})$, where $1\leq i\leq d$, in $k[x]$ are pairwise comaximal, the
Chinese Remainder Theorem yields
\[
R \cong \prod_{i=1}^d R_i, \where  R_i= \frac{k[x]}{(f_i(x)^{p^n})}.
\]
This implies that each $R$-module $M$ decomposes uniquely as $M = \bigoplus_{i=1}^dM_i$,
where $M_i$ is an $R_i$-module. Furthermore, it is easy to see that $R_i/(f_i(x)^s)$,
for $1\leq s\leq p^n$, is a complete list of indecomposable modules over $R_i$, and that
each $M_i$ has a unique decomposition into a direct sum of such modules. This is exactly
as predicted by the Krull--Schmidt theory.  The upshot is that we know `everything' about
the modules over the group algebras of finite cyclic groups.
\end{chunk}

All this is subsumed in the structure theory of modules over principal ideal rings.  By
the by, the finite cyclic groups are the source of group algebras of finite representation
type, in the following sense; see \cite[(4.4)]{Benson:1991a} for the appropriate references.

\begin{theorem}
\label{finite type}
If $k$ is an infinite field of characteristic $p$ and $G$ a finite group, then $\gring kG$
has finite representation type exactly when $G$ has cyclic Sylow $p$-subgroups. \hfill 
$\square$
 \end{theorem}

In some cases of infinite representation type, it is still possible to classify all the
indecomposable modules. The Klein group is one such. Let me give you a flavour of the
modules that arise over its group algebra. For the calculations, it is helpful to recall
a result on syzygies of indecomposable modules.

\begin{chunk}
\label{indecomposable:syzygies}
Let $(R,\fm,k)$ be a commutative artinian local ring and $E$ the injective hull of the
$R$-module $k$. Let $M$ be a finite $R$-module. Write $\Omega^1M$ for the first syzygy of $M$,
and $\Omega^{-1}M$ for the first co-syzygy of $M$.  These are defined by exact sequences
\[
0\to \Omega^1 M\to R^b\to M\to 0
  \also 0\to M\to E^c\to \Omega^{-1} M\to 0, \tag{$\dagger$}
\]
with $b=\rank_k(M/\fm M)$ and $c=\rank_k\soc M$.

The conclusion of the following exercise is valid for the syzygy module even when $R$ is
a Gorenstein ring of higher (Krull) dimension, as long as $M$ is also maximal
Cohen--Macaulay; this was first proved by J.~Herzog \citeyear{Her}.

\begin{Exercise}
Assume that $R$ is Gorenstein. Prove that when $M$ is indecomposable, so are $\Omega^1 M$
and $\Omega^{-1}M$.
\end{Exercise}
 
I cannot resist giving a sketch of the argument: Suppose $\Omega^1 M =U\oplus V$, with
$U$ and $V$ nonzero. Since $R$ is self-injective, neither $U$ nor $V$ can be free: if $U$
is free, then it is injective and hence splits from $R^b$ in the exact sequence
($\dagger$) above, and that cannot happen. Now, $\Hom R{-}R$ applied to ($\dagger$) yields
an exact sequence
 \[
0\to \dual M\to R^b \to \dual U\oplus \dual V \to 0.
\]
This presents $\dual M$ as the first syzygy of $\dual U\oplus \dual V$ (why?); that is,
\[
\dual M= \Omega^1(\dual U\oplus\dual  V) = \Omega^1(\dual U) \oplus \Omega^1(\dual V).
\]
Note that the modules $\Omega^1(\dual U)$ and $\Omega^1(\dual V)$ are nonzero: if
$\Omega^1(\dual U)=0$, then $\pdim_R(\dual U)$ is finite, so $\dual U$ is free,
and hence $U$ is free, a contradiction.  It follows that the same is true even
after we dualize them. Applying $\Hom R-R$ to the equality above gives us
\[
\dual{\dual M} = \dual{\Omega^1(\dual U)} \oplus \dual{\Omega^1(\dual V)}
\]
Since $M\cong \dual{(\dual M)}$, one obtains that $M$ is indecomposable.
 \end{chunk}

Now we turn to indecomposable modules over the Klein group.

\begin{chunk}{\obf The Klein group.}
\label{modules:klein}
Let $k$ be a field of characteristic $2$ and let $G$ be $\BZ_2\times \BZ_2$, the Klein
group.  Let $R$ denote its group algebra over $k$, so $R=k[y_1,y_2]/(y_2^2{-}1,y_2^2{-}1)$.

This $k$-algebra looks more familiar once we change variables: setting $x_i=y_i-1$ one
sees that $R=k[x_1,x_2]/(x_1^2,x_2^2)$; a local zero dimensional complete intersection
with maximal ideal $\fm=(x_1,x_2)$. Note that $R$ is Gorenstein, so $R\cong \Hom kRk$
and, for any $R$-module $M$, one has $\dual M\cong \Hom kMk$, where
$\dual{(-)}=\Hom RMR$.  I will use these remarks without ado.

For each positive integer $n$, let $M_n$ denote $\Omega^n(k)$, the $n$-th syzygy of $k$. I
claim that in the infinite family $\{\dots, M_2,M_1, k, \dual{(M_1)},\dual{(M_2)},\dots\}$ no
two modules are isomorphic and that each is indecomposable.

Indeed, a repeated application of Exercise \eqref{indecomposable:syzygies} yields that
each $M_n$ is indecomposable, and hence also that $\dual{(M_n)}$ is indecomposable,
since $\dual{\dual{(M_n)}}\cong M_n$. As to the remaining assertion: for $i=1,2$, let
$R_i=k[x_i]/(x_i^2)$. The minimal $R_i$-free resolution of $k$ is
\[
F_i=\quad \cdots \arto{x_i} R_i \arto{x_i} R_i \arto{x_i} R_i\to 0
\]
Since $R=R_1\otimes_kR_2$, the complex of $R$-modules $F_1\otimes_kF_2$ is the minimal
free resolution of the $R$-module $k$. It follows that the $n$-th Betti number of $k$ is $n+1$.
Thus, for any positive integer $n$, the $n$-th syzygy $M_n$ of $k$ is defined by an exact
sequence
\begin{equation*}
0\to M_n \to R^n \arto{\dd_{n-1}} R^{n-1}\to
        \cdots \to R^2 \arto{\dd_1} R\to  k\to 0,\tag{$\dagger$}
\end{equation*}
with $\dd_i(R^{i+1})\subseteq \fm R^i$ for each $i$.  It follows that $\rank_k M_n=2n+1$,
and hence also that $\rank_k\dual{(M_n)}=2n+1$.  Therefore, to settle the claim that the
modules in question are all distinct, it remains to verify that the $R$-modules $M_n$ and
$\dual{(M_n)}$ are not isomorphic.  These modules appear in exact sequences
\[
0\to M_n \to R^n \arto{\dd_{n-1}} R^{n-1} \qquad\text{and} \qquad
 0\to \dual{(M_n)}\to R^{n+1}\arto{\dual{\dd}_{n+1}} R^{n+2}.
\]
The one on the right is obtained from $$ R^{n+2}\arto{\dd_{n+1}} R^{n+1}\to M_n\to 0,$$
keeping in mind that $\dual{R}\cong R$.  Since $\dd_{n-1}(R^n)\subseteq \fm R^{n-1}$ and
$\dual{\dd}_{n+1}(R^{n+1})\subseteq \fm R^{n+2}$, the desired conclusion is a consequence
of:

\begin{Exercise}
  Let $(R,\fm,k)$ be a local ring. If $0\to K\to R^b\arto{f} R^c$ is an exact sequence of
  $R$-modules with $f(R^b)\subseteq \fm R^c$, then
\[
\soc K = \soc{R^b} = \soc R^b.
\]
\end{Exercise}

This completes the justification that the given family consists of nonisomorphic
indecomposables. In this process we found that $\rank_kM_n=2n+1=\rank_k\dual{(M_n)}$.  It
turns out that the $M_n$, their $k$-duals, and $k$ are the only indecomposables of odd rank;
here is a sketch of the proof. Exercise: fill in the details.

Let $M$ be an indecomposable $R$-module with $\rank_kM=2n+1$ for some integer $n$.  In
particular, $M\not\cong R$, and so Exercise \eqref{indecomposable:socleR} tells us that
$(xy)M=0$, so $\fm^2M=0$ and hence $\fm M\subseteq \soc M$; the opposite inclusion also
holds, by Exercise \eqref{indecomposable:socleM}, hence $\fm M=\soc M$. Thus, one has an
exact sequence of $R$-modules
\[
0\to \soc M \to M \to M/\fm M \to 0
\]
Now we use Exercise \eqref{indecomposable:syzygies}; in the notation there, from the exact
sequence above one deduces that either $b\leq n$ or $c\leq n$.  In the former case
$\rank_k(\Omega^1 M)\leq 2n-1$ and in the latter $\rank_k(\Omega^{-1}M)\leq 2n-1$.  In any
\vadjust{\goodbreak}%
case, the ranks of $\Omega^1 M$ and $\Omega^{-1}M$ are odd. Now an induction on rank
yields that $M$ belongs to the family of indecomposable $R$-modules that we have already
constructed.

At this point, we know all the indecomposable $R$-modules of odd rank. The ones of even
rank are harder to deal with. To get an idea of what goes on here, solve:

\begin{Exercise}
  Prove that every rank $2$ indecomposable $R$-module is isomorphic to a member of the
  family of cyclic $R$-modules
\[
V_{(\alpha_1,\alpha_2)} = \frac{R}{(\alpha_1 x_1 + \alpha_2 x_2,xy)}, \quad \text{where
   $(\alpha_1,\alpha_2)\ne (0,0)$}.
\]
Moreover, $V_{(\alpha_1,\alpha_2)}\cong V_{(\beta_1,\beta_2)}$ if and only if $(\alpha_1,\alpha_2)$
and$(\beta_1,\beta_2)$ are proportional.
\end{Exercise}
 
Thus, the nonisomorphic indecomposable $R$-modules of rank $2$ are parametrized by the
projective line over $k$; it turns out that this is the case in any even rank, at least when
$k$ is algebraically closed.  This classification of the indecomposable modules over the
Klein group goes back to Kronecker; see \cite{Al} or \cite[(4.3)]{Benson:1991a}
for a modern treatment.
\end{chunk}

This discussion shows that while the group algebra of $\BZ_2\times\BZ_2$ in
characteristic $2$ is not of finite type, in any given rank all but finitely many of its
indecomposable modules are contained in a one-parameter family. More generally, by
allowing for finitely many one-parameter families in each rank, one
obtains the notion of a
\textit{tame} algebra. Tame group algebras $\gring kG$ are completely classified: the
characteristic of $k$ is $2$, and the Sylow $2$-subgroups of $G$ are isomorphic to one of
the following groups: Klein, dihedral, semidihedral, or generalized quaternion. See
\cite[(4.4.4)]{Benson:1991a}. The significance of this result lies in
that every finite-dimensional $k$-algebra that is neither of finite
type nor tame is \textit{wild}, which 
implies that the set of isomorphism classes of its finite-rank indecomposable modules
contains representatives of the indecomposable modules over a tensor algebra in two variables.

\begin{ramble}
  There is a significant parallel between module theory over finite group algebras and
  over artinian commutative Gorenstein rings; see the discussion around Theorem
  \eqref{proj=inj}. In fact, this parallel extends to the category of maximal
  Cohen--Macaulay modules over commutative complete local Gorenstein rings. For example,
  analogous to Theorem \eqref{finite type}, among this class of rings those of
  \textit{finite Cohen--Macaulay type} (which means that there are only finitely many
  isomorphism classes of indecomposable maximal Cohen--Macaulay modules) have been
  completely classified, at least when the ring contains a field.  
A systematic exposition of this result can be found
  in \cite{Yo}.  The next order of complexity
  beyond finite Cohen--Macaulay type is bounded Cohen--Macaulay type, which is a topic of
  current research: see \cite{LW}.
\end{ramble}

The rest of this section describes a few basic constructions, like
tensor products and homomorphisms, involving modules over group algebras.

\begin{chunk}{\obf Conjugation.}
\label{conjugation}
Over a noncommutative ring, the category of left modules can be drastically different
from that of right modules. For example, there exist rings over which every left
module has a finite projective resolution, but not every right module does.  Thus, in
general, one has to be very careful vis-\`a-vis left and right module structures.

However, in the case of group algebras, each left module can be endowed with a natural
structure of a right module, and vice versa.  More precisely, if $M$ is a \textit{left}
$\gring kG$-module, then the $k$-vector space underlying $M$ may be viewed as a
\textit{right} $\gring kG$-module by setting
\[
m\cdot g = g^{-1}m \quad\text{for each $g\in G$ and $m\in M$}.
\]
For this reason, when dealing with modules over group algebras, one can afford to be lax
about whether they are left modules or right modules.  This also means, for instance, that
a left module is projective (or injective) if and only if the corresponding
right module has the same property.

This is similar to the situation over commutative rings: each left module $N$ over a
commutative ring $R$ is a right module with multiplication
\[
n\cdot r = rn \quad\text{for each $r\in R$ and $n\in N$}.
\]

There is an important distinction between the two situations: over $R$, the module $N$
becomes an $R$-bimodule with right module structure as above. However, over $\gring kG$,
the module $M$ with prescribed right module structure is not a bimodule.
 \end{chunk}

\begin{chunk}{\obf Tensor products.}
\label{tensors}
Over an arbitrary ring, one cannot define the tensor product of two left modules. However,
if $M$ and $N$ are two left modules over a group algebra $\gring kG$, one can view
$M$ as a right module via conjugation \eqref{conjugation} and make sense of
$M\otimes_{\gring kG}N$. But then this tensor product is \textit{not} a $\gring kG$-module,
because $M$ and $N$ are not bimodules.  In this respect, the group ring
behaves like any old ring.

There is another tensor product construction, a lot more important when dealing with group
algebras than the one above, that gives us back a $\gring kG$-module. To describe
it, we return briefly to the world of arbitrary $k$-algebras.

Let $R$ and $S$ be $k$-algebras and let $M$ and $N$ be (left) modules over $R$ and $S$,
respectively. There is a natural left ($\tensor RkS$)-module structure on $\tensor MkN$
with
\[
(\tensor rks)\cdot(\tensor mkn)=\tensor{rm}k{sn}.
\]

Now let $M$ and $N$ be left $\gring kG$-modules. The preceding recipe provides an action
of $\tensor {\gring kG}k{\gring kG}$ on $\tensor MkN$.  This restricts, via the diagonal
map \eqref{gring:diagonal}, to a left $\gring kG$-module structure on $\tensor MkN$. Going
through the definitions one finds that
\[
g\cdot (\tensor mkn) = \tensor {gm}k{gn},
\]
for all $g\in G$, $m\in M$ and $n\in N$. It is worth remarking that the `twisting' map
\begin{align*}
M&\otimes_kN \arto{\cong} N\otimes_kM,\\
(m&\otimes_kn) \mapsto (n\otimes_km),
\end{align*}
which is bijective, is $\gring kG$-linear.  
 \end{chunk}

\begin{ramble}
  To a commutative algebraist, the tensor product $M\otimes_kN$ has an unsettling feature:
  it is taken over $k$, rather than over $\gring kG$. However, bear in mind that the
  $\gring kG$-module structure on $M\otimes_kN$ uses the diagonal homomorphism. The
  other possibilities, namely acquiring the structure from $M$ or from~$N$, don't give us
  anything nearly as useful. For instance, $M\otimes_kN$ viewed as a $\gring kG$-module
  via its left-hand factor is just a direct sum of copies of $M$.
\end{ramble}

\begin{chunk}{\obf Homomorphisms.}
\label{homs}
Let $M$ and $N$ be left $\gring kG$-modules. One can then consider $\Hom{\gring kG}MN$,
the $k$-vector space of $\gring kG$-linear maps from $M$ to $N$. Like the tensor product
over $\gring kG$, this is not, in general, a $\gring kG$-module. Note that since the
$\gring kG$-module $k$ is cyclic with annihilator $\aug G$, and $\aug G$ is generated as
an ideal by elements $g-1$, one has
\[
\Hom{\gring kG}kM = \{m\in M\mid gm=m\}.
\]
The $k$-subspace on the right is of course $M^G$, the set of $G$-invariant elements in~$M$.

As with $M\otimes_kN$, one can endow the $k$-vector space $\Hom kMN$ with a canonical left
$\gring kG$-structure. This is given by the following prescription: for each $g\in G$,
$\alpha\in \Hom kMN$, and $m\in M$, one has
\[
(g\cdot \alpha)(m) = g\alpha(g^{-1} m).
\]
In particular, $g\cdot\alpha=\alpha$ if and only if $\alpha(gm)=g\alpha(m)$; that is to
say,
\[
\Hom{\gring kG}MN = {\Hom kMN}^G.
\]
Thus the homomorphisms functor $\Hom{\gring kG}MN$ is recovered as the $k$-subspace of
$G$-invariant elements in $\Hom kMN$. This identification leads to the following
Hom-Tensor adjunction formula:
\[
\Hom{\gring kG}{\tensor LkM}N \cong \Hom {\gring kG}L{\Hom kMN}.
\]
This avatar of Hom-Tensor adjunction is very useful in the study of modules over group
algebras; see, for example, the proof of \eqref{proj:ideal}.
\end{chunk}

\begin{ramble}
Let $G$ be a finite group such that the characteristic of $k$ is coprime to $|G|$, and
let $0\to L\to M\to N\to 0$ be an exact sequence $\gring kG$-modules.
Applying $\Hom{\gring kG}k-$ to it yields, in view of Maschke's theorem \eqref{maschke:real},
an exact sequence
\[
0\to L^G \to M^G \to N^G\to 0.
\]
This is why invariant theory in characteristics coprime to $|G|$ is so
drastically different
from that in the case where the characteristic of $k$ divides $|G|$.
\end{ramble}

\begin{chunk}{\obf\hskip-1pt A technical point.}
Let $M$ be a left $\gring kG$-module and set $\thickmuskip5mu minus 5mu\dual M=\Hom kMk$.  One has two choices for
\label{standard isos}
a left $\gring kG$-module structure on $\dual M$: one given by specializing the discussion
in \eqref{homs} to the case where $N=k$, and the other by conjugation\emdash see
\eqref{conjugation}\emdash from the natural \textit{right} module structure on $\dual M$.  A
direct calculation reveals that they coincide. What is more, these modules have the
property that the canonical maps of $k$-vector spaces
\[
\begin{gathered}
   \begin{aligned}
      M& \to \dual {\dual M}      \\
      m&\mapsto \big(f\mapsto f(m)\big) 
    \end{aligned}\qquad\qquad
 \begin{aligned}
 N&\otimes_k \dual M \to \Hom kMN\\
 n&\otimes_kf \mapsto \big(m\mapsto f(m)n\big)
   \end{aligned}
\end{gathered}
\]
are $\gring kG$-linear.  These maps are bijective when $\rank_kM$ is finite.
\end{chunk}

\begin{ramble}
  Most of what I said from \eqref{conjugation} onward applies, with appropriate
  modifications, to arbitrary Hopf algebras. For example, given modules $M$ and $N$ over a
  Hopf algebra $H$, the tensor product $M\otimes_kN$ is also an $H$-module with
\[
h\cdot(m\otimes_kn) =\sum_{i=1}^n h_i'm\otimes_kh_i''n ,
  \quad\text{where}\ \Delta(h) = \sum_{i=1}^n h_i'\otimes_kh_i'' .
\]
There are exceptions; for example, over a group algebra $M\otimes_kN\cong N\otimes_kM$;
see \eqref{tensors}.  This holds over $H$ only when $\sum_{i=1}^n h_i'\otimes_kh_i''
=\sum_{i=1}^n h_i''\otimes_kh_i'$, that is to say, when the diagram
\[
\xymatrixrowsep{1pc}
\xymatrixcolsep{1.5pc} 
\xymatrix{
     &H\ar@{->}[ld]_{\Delta} \ar@{->}[rd]^{\Delta} & \\
H\otimes_k H \ar@{->}[rr]_{\tau}& &
  H\otimes_k H }
\]
commutes,
where $\tau(h'\otimes_kh'')=(h''\otimes_kh')$. Such an $H$ is said to be
\textit{cocommutative}.
\end{ramble}
\section{Projective Modules}
The section focuses on projective modules over group algebras. First, I address
the question: When is every module over the group algebra projective? In other words, when
is the group algebra \textit{semisimple}? Here is a complete answer, at least in the case
of a finite group.

\begin{theorem}
\label{maschke}
Let $G$ be a finite group. The following conditions are equivalent:
\begin{enumerate}
\item[{\rm(i)}] The group ring $\gring kG$ is semisimple.
\item[{\rm(ii)}] $k$, viewed as a $\gring kG$-module via the augmentation, is projective.
\item[{\rm(iii)}] The characteristic of $k$ is coprime to $|G|$.
\end{enumerate}
\end{theorem}
\begin{proof}
(i) $\implies$ (ii) is a tautology.
  
\smallskip\noindent
(ii) $\implies$ (iii): As $k$ is projective, the augmentation homomorphism $\eps\col
\gring kG\to k$, being a surjection, has a $\gring kG$-linear section $\sigma\col k\to
\gring kG$. Write $\sigma(1) = \sum_{g\in G}a_gg$, with $a_g$ in $k$.  Fix an element
$h\in G$. Note that $\sigma(1)=\sigma(h\cdot 1)=h\cdot\sigma(1)$, where the first equality
holds because $\gring kG$ acts on $k$ via $\eps$, the second by the $\gring
kG$-linearity of $\sigma$. This explains the first equality below:
\[
\sum_{g\in G}a_gg = \sum_{g\in G}a_g(hg) = \sum_{g\in G}a_{h^{-1}g}g.
\]
The second is just a reindexing. The elements of $G$ are a basis for the group algebra, so
the equality above entails $a_{h^{-1}}=a_1$. This holds for each $h\in G$, so
\[
1=\eps(\sigma(1)) = a_1\sum_{g\in G}\eps(g) = a_1\sum_{g\in G}1 = a_1|G|.
\]
In particular, the characteristic of $k$ is coprime to $|G|$.

\smallskip\noindent
(iii) $\implies$ (i): Let $M$ be a $\gring kG$-module, and pick a surjection $P\tra M$ with
$P$ projective.  Maschke's theorem \eqref{maschke:real} provides that every short exact
sequence of $\gring kG$-modules splits; equivalently, that every surjective homomorphism
is split. In particular, $P\tra M$ splits, so $M$ is a direct summand of $P$, and
hence projective.
  \end{proof}

\begin{Exercise}
  A commutative ring is semisimple if and only if it is a product of fields.
\end{Exercise}
 
The last result dealt with modules en masse; now the focus is shifted to individual
modules. 

\subsection*{Stability properties of projective modules.}
The gist of the following paragraphs is that many of the standard functors of interest
preserve projectivity. A crucial, and remarkable, result in this direction is 

\begin{theorem}\label{proj:ideal}
  Let $G$ be a group and $P$ a projective $\gring kG$-module. For any $\gring kG$-module
  $X$, the $\gring kG$-modules $P\otimes_kX$ and $X\otimes_kP$ are projective.
\end{theorem}

Take note that the tensor product is over $k$, as it must be, for such a conclusion is
utterly wrong were it over $\gring kG$. This theorem underscores the point raised
in \eqref{tensors} about the importance of this tensor product in the module theory of
group algebras; the other results in this section are all formal consequences of this one.

\begin{ramble}
  There is another way to think about Theorem \eqref{proj:ideal}: one may view the entire
  category of $\gring kG$-modules as a `ring' with direct sum and tensor product over $k$
  playing the role of addition and multiplication respectively; the unit is $k$, and the
  commutativity of the tensor product means that this is even a
  `commutative' ring. (With suitable compatibility conditions, such data define a \textit{symmetric monoidal
      category}.)  In this language, the theorem above is equivalent
  to the statement 
  that the subcategory of projective modules is an ideal.
 \end{ramble}

\begin{proof}[Proof of Theorem \eqref{proj:ideal}]
  I will prove that $P\otimes_kX$ is projective. A similar argument works for
  $X\otimes_kP$; alternatively, note that it is isomorphic to $P\otimes_kX$, by
  \eqref{tensors}.
  
  One way to deduce that $P\otimes_kX$ is projective is to invoke the following
  isomorphism from \eqref{homs}, which is natural on the category of left $\gring
  kG$-modules:
\[
\Hom{\gring kG}{\tensor PkX}- \cong \Hom {\gring kG}P{\Hom kX-}.
\]

Perhaps the following proof is more illuminating: by standard arguments one reduces to the
case where $P=\gring kG$.  Write $X\nat$ for the $k$-vector space underlying $X$.  Now, by
general principles, the inclusion of $k$-vector spaces $X\nat \subset \gring kG\otimes_kX$,
defined by $x\mapsto 1\otimes_kx$, induces a $\gring kG$-linear map
\[
\gring kG\otimes_k X\nat \to \gring kG\otimes_kX , \where g\otimes_kx\mapsto
g(1\otimes_kx) = g\otimes_k{gx}.
\]
The action of $\gring kG$ on $\gring kG\otimes_kX\nat$ is \textit{via the left-hand factor}.
An elementary calculation verifies that the map below, which is $\gring kG$-linear, is its
inverse:
\[
\gring kG\otimes_kX \to \gring kG\otimes_k X\nat, \where g\otimes_kx\mapsto
g\otimes_k(g^{-1}x).
\]
Therefore, the $\gring kG$-modules $\gring kG\otimes_kX$ and $\gring kG\otimes_k{X\nat}$
are isomorphic. It remains to note that the latter module is a direct sum of copies of
$\gring kG$.
 \end{proof}
 
 One corollary of Theorem \eqref{proj:ideal} is the following recognition principle
 for semisimplicity of the group algebra; it extends to arbitrary groups the equivalence
 of conditions (i) and (ii) in Theorem \eqref{maschke}.

\begin{lemma}
\label{semisimple}
Let $G$ be a group. The following conditions are equivalent.
\begin{enumerate}
\item[{\rm(i)}] $\gring kG$ is semisimple;
\item[{\rm(ii)}] the $\gring kG$-module $k$ is projective.
\end{enumerate}
\end{lemma}
\begin{proof}
  The nontrivial implication is that (ii) $\implies$ (i). As to that, it follows from
  Theorem \eqref{proj:ideal} that $k\otimes_kM$ is projective for each $\gring kG$-module
  $M$, so it remains to check that the canonical isomorphism $k\otimes_kM \to M$ is
  $\gring kG$-linear. Note that this is something that needs checking for
  the $\gring kG$-action on $k\otimes_kM$ is via the diagonal; see \eqref{tensors}.
\end{proof}

\begin{ramble}
  Lemma \eqref{semisimple}, although not its proof, is reminiscent of a phenomenon 
 encountered in the theory of commutative local rings: Over such a ring, the residue
  field is often a `test' module.  The Auslander--Buchsbaum--Serre characterization of
  regularity is no doubt the most celebrated example.  It says that a
 noetherian commutative local ring $R$, with residue field $k$, 
is regular if and only if
the $R$-module $k$ has finite projective dimension.

There are analogous results that characterize complete intersections
(Avramov and Gulliksen) and Gorenstein rings (Auslander and Bridger).
  
There is however an important distinction between a group algebra over $k$ and a local
ring with residue field $k$: over the latter, $k$ is the only simple module, whilst the
former can have many others.  From this perspective, Lemma \eqref{semisimple} is rather
surprising. The point is that an arbitrary finite-dimensional algebra is semisimple if and
only if {\em every} simple module is projective; the nontrivial implication holds because
each finite module has a composition series.
 \end{ramble}

\begin{theorem}\label{proj:functors}
  Let $G$ be a finite group. For each finite $\gring kG$-module $M$, the following $\gring
  kG$-modules are projective  simultaneously: $M$, $M\otimes_kM$, $\dual M\otimes_kM$,
  $M\otimes_k\dual M$, $\Hom kMM$, and $\dual M$.
 \end{theorem}
\begin{proof}
  It suffices to verify: $M$, $M\otimes_kM$, and $\dual M\otimes_kM$ are simultaneously
  projective.
  
  Indeed, applied to $\dual M$ that would imply, in particular, that $\dual M$ and
  $\dual{(\dual M)}\otimes_k\dual M$ are simultaneously projective.  Now, $\dual{(\dual
    M)}\cong M$, since $\rank_kM$ is finite, and $M\otimes_k\dual M\cong \dual M\otimes_k
  M$, by the discussion in \eqref{tensors}. Thus, one obtains the simultaneous projectivity
  of all the modules in question, except for $\Hom kMM$. However, the finiteness of
  $\rank_kM$ implies this last module is isomorphic to $M\otimes_k\dual M$.
  
  As to the desired simultaneous projectivity, it is justified by the diagram
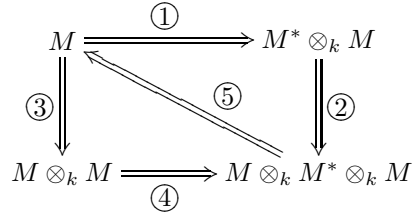
\[
\xymatrixrowsep{3pc} 
\xymatrixcolsep{3pc} 
\xymatrix{ 
M \ar@{=>}[r]^(.4){\xy*+{1}*\cir<5pt>{}\endxy} \ar@{=>}[d]_{\xy*+{3}*\cir<5pt>{}\endxy}
  & \dual M \otimes_k M \ar@{=>}[d]^{\xy*+{2}*\cir<5pt>{}\endxy} \\
  M\otimes_kM\ar@{=>}[r]_(.4){\xy*+{4}*\cir<5pt>{}\endxy}
  &\ar@{=>}[ul]_(.4){\xy*+{5}*\cir<5pt>{}\endxy} M\otimes_k\dual M\otimes_kM}
\]
 where $X\Rightarrow Y$ should be read as `if $X$ is projective, then so is $Y$'.
Implications (1)--(4) hold by Theorem \eqref{proj:ideal}.  As to (5), the
natural maps of $k$-vector spaces
\begin{alignat*}2
M &\to \Hom kMM\otimes_kM &&\to M \\
m &\mapsto 1\otimes_k m\ \text{and}\ \alpha\otimes_k m &&\mapsto \alpha(m)
\end{alignat*}
are $\gring kG$-linear, and exhibit $M$ as a direct summand of $\tensor{\Hom kMM}kM$.
However, as remarked before, the $\gring kG$-modules $\Hom kMM$ and $M\otimes_k\dual M$
are isomorphic, so $M$ is a direct summand of $M\otimes_k\dual M\otimes_kM$.
\end{proof}

\subsection*{Projective versus Injectives.}

So far, I have focused on projective modules, without saying anything at all about
injective, or flat, modules. Now, a commutative algebraist well knows that projective
modules and injective modules are very different beasts. There is, however, one exception.

\begin{exercise}
  Let $R$ be a commutative noetherian local ring. Prove that when $R$ is zero-dimensional and
  Gorenstein, an $R$-module is projective if and only if it is injective. Conversely, if
  there is a nonzero $R$-module that is both projective and injective, then $R$ is
  zero-dimensional and Gorenstein.
\end{exercise}

The preceding exercise should be compared with the next two results.

\begin{theorem}
\label{proj=inj}
Let $G$ be a finite group and $M$ a finite $\gring kG$-module. The following conditions
are equivalent:
\begin{enumerate}
\item[{\rm(i)}] $M$ is projective;
\item[{\rm(ii)}] the flat dimension of $M$ is finite;
\item[{\rm(iii)}] $M$ is injective;
\item[{\rm(iv)}] the injective dimension of $M$ is finite. 
 \end{enumerate}
 \end{theorem}

These equivalences hold for any $\gring kG$-module, finite or
 not; see \cite{Be:fm}. 

The preceding theorem has an important corollary.

\begin{corollary}\label{gring:gor}
The group algebra of a finite group is self-injective.\hfill$\square$
 \end{corollary}
 
 There are many other proofs, long and short, of this corollary; see
 \cite[(3.1.2)]{Benson:1991a}.  Moreover, it is an easy exercise (do
 it) to deduce Theorem \eqref{proj=inj} from it.{\looseness=-1\par}

\begin{ramble}
  Let $G$ be a finite group. Thus, the group algebra $\gring kG$ is finite-dimensional
  and, by the preceding corollary, injective as a module over itself.  These properties
  may tempt us commutative algebraists to proclaim: $\gring kG$ is a zero-dimensional
  Gorenstein ring.  And, for many purposes, this is a useful point of view, since module
  theory over a group algebra resembles that over a Gorenstein ring; Theorem
  \eqref{proj=inj} is one manifestation of this phenomenon. By the by, there are diverse
  extensions of the Gorenstein property for commutative rings to the noncommutative 
  setting: Frobenius rings, quasi-Frobenius rings, symmetric rings, self-injective rings,
  etc.
\end{ramble}

The proof of Theorem \eqref{proj=inj} is based on Theorem \eqref{proj:functors} and an
elementary observation about modules over finite-dimensional algebras.

\begin{lemma}
  Let $R$ be a $k$-algebra with $\rank_kR$ finite. For each finite left $R$-module $M$,
  one has $\pdim_RM = \fdim_R M = \injdim_{R^{\rm{op}}}\dual M$. 
\end{lemma}

\begin{proof}
Since $\rank_k M$ is finite, $\dual{(\dual M)}\cong M$, so it suffices to prove the
equivalence of the conditions

\begin{enumerate}
\item[{\rm(i)}] $M$ is projective;
\item[{\rm(ii)}] $M$ is flat;
\item[{\rm(iii)}] the right $R$-module $\dual M$ is injective.
 \end{enumerate}  
 The implication (i) $\implies$ (ii) is immediate and hold for
 all rings.  The equivalence
 (ii) $\iff$ (iii) is a consequence of the standard adjunction isomorphism
\[
\Hom k{-\otimes_RM}k\cong \Hom R-{\dual M}
\]
and is valid for arbitrary $k$-algebras. 

(iii) $\implies$ (i): Since $M$ is finite over $R$, one can construct a surjective map
$\pi\col R^n\tra M$. Dualizing this yields an inclusion $\dual\pi\col \dual M
\hookrightarrow \dual{(R^n)}$ of right $R$-modules.  This map is split because $\dual M$
is injective, and hence $\dual{\dual\pi}$ is split.  Since $\rank_kR$ and $\rank_kM$ are
both finite, $\dual{\dual \pi}=\pi$, so that $\pi$ is split as well. Thus, $M$ is
projective, as claimed.
 \end{proof}

\begin{proof}[Proof of Theorem \eqref{proj=inj}]
  Theorem \eqref{proj:functors} yields that $M$ is projective if and only if $\dual M$ is
  projective, while the lemma above implies that $\dual M$ is projective if and only if
  $\dual{(\dual M)}$ is injective, i.e.,~$M$ is injective. This settles (i) $\iff$ (iii).
  
  That (i) $\implies$ (ii) needs no comment. The lemma above contains (ii)
  $\iff$ (iv); moreover, it implies that to verify (ii) $\implies$ (i), one may assume
  $\pdim_RM$ finite, that is to say, there is an exact sequence
\[
0\to P_n \arto{\dd_n} P_{n-1}\arto{\dd_{n-1}}\cdots \to P_0 \to M\to 0,
\]
where each $P_i$ is finite and projective; see \eqref{resolutions}. If $n\geq 1$, then,
since $P_n$ is injective by the already verified implication (i) $\implies$ (iii), the
homomorphism $\dd_n$ splits, and one obtains an exact sequence
\[
0\to \dd_{n-1}(P_{n-1}) \to P_{n-2}\to \cdots \to P_0\to M\to 0.
\]
In this sequence $\dd_{n-1}(P_{n-1})$, being a direct summand of $P_{n-1}$, is projective,
and hence injective. An iteration of the preceding argument yields that $M$ is a direct
summand of $P_0$, and hence projective. 
\end{proof}

\begin{ramble}
  The small finitistic left global dimension of a ring $R$ is defined as
\[
\sup\,\{\pdim_RM \mid \text{$M$ a finite left $R$-module with $\pdim_RM<\infty$}.\}
\]
One way of rephrasing Theorem \eqref{proj=inj} is to say that this number is zero when $R$
is a finite group algebra.  Exercise: Prove that a similar result holds also for modules
over commutative artinian rings.  However, over arbitrary finite-dimensional algebras, the
small finitistic global dimension can be any nonnegative integer.  A conjecture of Bass
\citeyear{Ba} and Jans \citeyear{Ja}, which remains open, asserts that this number is finite; look
up \cite{Ha} for more information on this topic.
 \end{ramble}
 
\subsection*{Hopf algebras.}
Theorem \eqref{proj:ideal} holds also for modules over any finite-dimen\-sional Hopf
algebra; the proof via the adjunction isomorphism does not work, but the other one does.
However, I found it a nontrivial task to pin down the details, and I can recommend
it to you as a good way to gain familiarity with Hopf algebras. Given this, it is not hard
to see that for \textit{cocommutative} Hopf algebras, the analogues of theorems
\eqref{proj:functors} and \eqref{proj=inj}, and Corollary \eqref{gring:gor}, all hold;
the cocommutativity comes in because in the proof of \eqref{proj:functors} I used the 
fact that tensor products are symmetric; confer with the discussion in \eqref{standard isos}.
 
\section{Structure of Projectives}
 
So far, I have not addressed the natural question: what are the projective modules over
the group algebra? In this section, I tabulate some crucial facts concerning
these. Most are valid for arbitrary finite-dimensional algebras and are easier to state in
that generality; \cite{Al} is an excellent reference for this circle of ideas.

\medskip

\begin{chunk}{\obf Projective covers.}
\label{covers}
Let $R$ be a ring and $M$ a finite $R$-module.  A \textit{projective cover} of $M$ is a
surjective homomorphism $\pi\col P\to M$ with $P$ a projective $R$-module and such that each
homomorphism $\sigma\col P\to P$ that fits in a commutative diagram
\[
\xymatrixrowsep{1.5pc} 
\xymatrixcolsep{2pc} 
\xymatrix{
  P \ar[dr]^{\pi} \ar[rr]^{\sigma}
&    & P\ar[dl]_{\pi} \\
  &M & }
\]
is bijective, and hence an automorphism.  It is clear that projective covers, when they
exist, are unique up to isomorphism. Thus, one speaks of \textit{the} projective cover of
$M$.  Often $P$, rather than $\pi$, is thought as being the projective cover of
$M$, although this is an abuse of terminology.

Among surjective homomorphisms $\kappa\col Q\to M$ with $Q$ a projective $R$-module,
projective covers can be characterized by either of the properties:
\begin{enumerate}
\item [{\rm(i)}]
$Q/JQ\cong M/JM$, where $J$ is the Jacobson radical of $R$;
\item[{\rm(ii)}]
$Q$ is minimal with respect to direct sum decompositions.
\end{enumerate}

When $R$ is a noetherian ring over which every finite $R$-module has a projective cover, it
is easy to see that a projective resolution 
\[
\bsP\col \cdots \to P_n \arto{\dd_n}P_{n-1}\arto{\dd_{n-1}}\cdots \arto{\dd_1}P_0\to 0
\]
of $M$ so constructed that $P_n$ is a projective cover of $\Ker(\dd_{n-1})$ is unique up
to isomorphism of complexes of $R$-modules. Such a $\bsP$ is called the \textit{minimal
  projective resolution} of $M$. Following conditions (i) and (ii) above, the minimality
can also be characterized by either the property that $\dd(\bsP) \subseteq J\bsP$, or that
$\bsP$ splits off from any projective resolution of $M$.

Projective covers exist for each finite $M$ in two cases of interest: when $R$ is a finite-dimensional $k$-algebra, and when $R$ is a (commutative) local ring. This is why these two
classes of rings have a parallel theory of minimal resolutions.
\end{chunk}
 
\begin{chunk}{\obf Simple modules.}
\label{blocks}  
Let $R$ be a finite-dimensional $k$-algebra with Jacobson radical $J$, and let ${\mathcal
  P}$ and ${\mathcal S}$ be the isomorphism classes of indecomposable projective
$R$-modules and of simple $R$-modules, respectively.
\begin{enumerate}
\item[{\rm(a)}]
The Krull--Schimdt property holds for $R$, so every $P$ in $\mathcal P$ occurs as a
direct summand of $R$, and there is a unique decomposition
\[
\begin{gathered}
{\qquad}
{R \cong \bigoplus_{P\in \mathcal P} P^{e_R(P)}},
{\quad} {\text{with \,$e_R(P)\geq 1$.}}
\end{gathered}
\]
In particular, $R$ has only finitely many indecomposable projective modules.
\item [{\rm(b)}]
The simple $R$-modules are precisely the indecomposable modules of the semisimple
 ring $\wt R= R/J$ (verify this) so property (a) specialized to $\wt R$ reads
\[
\begin{gathered}
{\qquad}
{\wt R \cong \bigoplus_{S\in \mathcal S} S^{e_{\wt R}(S)}},
{\quad} {\text{with \,$e_{\wt R}(S)\geq 1$.}}
\end{gathered}
\]
\item[{\rm(c)}] The ring $\wt R$ in (b), being semisimple, is a direct sum of matrix rings over
  finite-dimensional division algebras over $k$; see \cite[XVII]{La}. Moreover, when $k$ is
  algebraically closed, these division algebras coincide with $k$ (why?), and we obtain
  that $e_{\wt R}(S) = \rank_kS$ for each $S\in\mathcal S$.
\item[{\rm(d)}] From (a)--(c) one obtains that the assignment $P\mapsto P/JP$ is a bijection between
  $\mathcal P$ and $\mathcal S$; in other words, there are as many indecomposable
  projective $R$-modules as there are simple $R$-modules. Moreover, $e_R(P)=e_{\wt
    R}(P/JP)$. When $k$ is algebraically closed, combining the last equality with that in
  (c) and the decomposition in (a) yields
\[
\rank_kR = \sum_{P\in\mathcal P} \rank_k(P/JP)\rank_k P.
\]
\end{enumerate}
\end{chunk}

I will illustrate the preceding remarks by describing the indecomposable projective
modules over certain finite group algebras.

\begin{chunk}{\obf Cyclic groups.}
\label{projectives:cyclic}  
This example builds on the description in \eqref{modules:cyclic} of modules over the group
algebra of a finite cyclic group $G$. We saw there that
\[
\gring kG \cong \prod_{i=1}^d \frac{k[x]}{( f_i(x)^{p^n})}.
\]
This is the decomposition that for general finite-dimensional algebras is a consequence of
the Krull--Schmidt property; see (\ref{blocks}.a).  For each $1\leq i\leq d$, set
$P_i=k[x]/(f_i(x)^{p^n})$. These $\gring kG$-modules are all projective, being
summands of $\gring kG$, indecomposable (why?), and no two of them are isomorphic (count
ranks, or look at their annihilators).  Moreover, as a consequence of the decomposition
above, any projective $\gring kG$-module is a direct sum of the $P_i$. Thus, there are
exactly $d$ distinct isomorphism classes of indecomposable projective $R$-modules.

Over any commutative ring, the only simple modules are the residue fields. Thus, the
simple modules over $\gring kG$ are $k[x]/(f_i(x))$ where $1\leq i\leq d$; in particular, there
are as many as there are indecomposable projectives, exactly as (\ref{blocks}.d) predicts.
\end{chunk}

Now I will describe the situation over finite abelian groups. Most of what I have to
say can be deduced from:

\begin{lemma}
\label{projectives:product}
Let $R$ and $S$ be finite-dimensional $k$-algebras, and set $T=R\otimes_kS$. 
Let $M$ and $N$ be $R$-modules. If $S$ is local with residue ring is $k$, and the induced
map $k\to S\to k$ is the identity, then
\begin{enumerate}
\item[{\rm(a)}] $M\cong N$ as $R$-modules if and only if $M\otimes_kS\cong N\otimes_kS$ as 
  $T$-modules;
\item[{\rm(b)}] the $R$-module $M$ is indecomposable if and only if the $T$-module $M\otimes_kS$ is;
\item[{\rm(c)}] $M$ is projective if and only if the $T$-module $M\otimes_kS$ is projective.
\end{enumerate}
In particular, the map $P\mapsto P\otimes_kS$ induces a bijection between the isomorphism
classes of indecomposable projective modules over $R$ and over $T$.
\end{lemma}

\begin{proof}
  To begin with, note that $M\otimes_kS$ and $N\otimes_kS$ are both left $R$-modules and
  also right $S$-modules, with the obvious actions.  Moreover, because of our hypothesis
  that the residue ring of $S$ is $k$, one has isomorphisms of $R$-modules
\[
M\cong (M\otimes_k S)\otimes_Sk \qquad\text{and}\qquad
N\cong (N\otimes_k S)\otimes_Sk.
\]
Now, the nontrivial implication in (a) and in (c)\emdash the one concerning descent\emdash is
settled by applying $-\otimes_Sk$. As to (b), the moot point is the ascent, so assume the
$R$-module $M$ is indecomposable and that $M\otimes_k S\cong U\oplus V$ as $T$-modules.
Applying $-\otimes_Sk$, one obtains isomorphisms of $R$-modules
\[
M\cong (M\otimes_k S)\otimes_Sk \cong (U\otimes_Sk) \oplus (V\otimes_Sk)
\]
Since $M$ is indecomposable, one of $U\otimes_Sk$ or $V\otimes_Sk$ is zero; say,
$U\otimes_Sk$ is $0$, that is to say, $U=U\fn$, where $\fn$ is the maximal ideal of $S$.
This implies $U=0$, because, $S$ being local and finite-dimensional over
$k$, the ideal $\fn$ is nilpotent.
\end{proof}

\begin{chunk}{\obf Finite abelian groups.}
\label{projectives:abelian}   
Again, we adopt that convention that $p$ is the characteristic of $k$ when the latter is
positive, and $1$ otherwise.

Let $G$ a finite abelian group, and write $|G|$ as $p^nq$, where $n$ is a nonnegative
integer and $p$ and $q$ are coprime.  Via the fundamental theorem on finitely generated
abelian groups this decomposition of $|G|$ translates into one of groups:
$G=A\oplus B$, where $A$ and $B$ are abelian, $|A|=p^n$, and $|B|=q$. Hence, $\gring
kG\cong \gring kA\otimes_k\gring kB$.

Now, $A\cong \BZ/(p^{e_1}\BZ)\oplus \cdots \oplus\BZ/(p^{e_m}\BZ)$, for nonnegative
integers $e_1,\dots,e_m$, so 
\[
\gring kA\cong \frac{k[y_1,\dots,y_m]}{(y_1^{p^{e_1}}{-}1,\dots,y_m^{p^{e_m}}{-}1)}
\]
The binomial theorem in characteristic $p$ yields $y_i^{p^{e_i}}-1 =
(y_i-1)^{p^{e_i}}$ for each $i$. Thus, it is clear that $\gring kA$ is an artinian local
ring with residue field $k$.

In the light of this and Lemma \eqref{projectives:product},
to find the indecomposable projectives over $\gring kG$, it suffices to find those
over $\gring kB$. 

When $B$ is cyclic, this information is contained in \eqref{projectives:cyclic}.  The
general case is more delicate. First, since $|B|$ is coprime to $p$, every
$\gring kG$-module is projective, so the indecomposables among them are precisely the
simple $\gring kB$-modules; see Theorem \eqref{maschke}. Now, as noted before, 
over any commutative ring the only simple modules are the residue fields. Thus, the
problem is to find the maximal ideals of $\gring kB$. Writing $B$ as $\BZ/(q_1\BZ)
\oplus \cdots \oplus\BZ/(q_n\BZ)$, one has
\[
\gring kB\cong \frac{k[x_1,\dots,x_n]}{(x_1^{q_1}{-}1,\dots,x_m^{q_n}{-}1)}.
\]
If $k$ is algebraically closed, there are $q_1\cdots q_n$ distinct maximal ideals,
and hence as many distinct indecomposable projectives. The general situation is
trickier.

By the by, if you use the method outlined above for constructing projective modules over a
cyclic group, the outcome will appear to differ from that given by
\eqref{projectives:cyclic}. Exercise: Reconcile them.
 \end{chunk}

\begin{chunk}{\obf $p$-groups.}
\label{pgroups}
As always, free $\gring kG$-modules are projective. When the characteristic of $k$ is $p$
and $G$ is a $p$-group, these are the only projectives over $\gring kG$.  This is thus
akin to the situation over commutative local rings, and the proof over this latter class
of rings given in \cite{Matsumura:1986a} carries over; the key ingredient is that, as
noted in \eqref{gring:pgroups}, the group algebra of a $p$-group is an artinian local
ring.
\end{chunk}

In general, the structure of projective modules over the group algebra is a lot more
complicated. However, the triviality of the projectives in the case of $p$-groups also has
implications for the possible ranks of indecomposable projectives over the group algebra
of an arbitrary group $G$.

\begin{chunk}{\obf Sylow subgroups.}
\label{projectives:sylows}
Let $p^d$ be the order of a $p$-Sylow subgroup of $G$. If a finite $\gring kG$-module $P$
is projective, then $p^d$ divides $\rank_kP$.

Indeed, for each $p$-Sylow subgroup $H\subseteq G$, the restriction of $P$ to the subring
$\gring kH$ of $\gring kG$ is a projective module, and hence a free module.  Thus, by the
preceding remark, $\rank_kP$ is divisible by $\rank_k\, \gring kH$, that is to say, by
$|H|$.
\end{chunk}

The numerological restrictions in \eqref{blocks} and \eqref{projectives:sylows} can be
very handy when hunting for projective modules over finite group algebras.
Here is a demonstration.

\begin{chunk}{\obf Symmetric group on three letters.}
 \label{projectives:example}
The symmetric group on three letters,
 $\Sigma_3$, is generated by elements $a$ and $b$, subject to the relations
\[
a^2=1, \quad b^3=1, \qnd ba=ab^2.
\]
Thus, $\Sigma_3=\{1,b,b^2,a,ab,ab^2\}$. It has three 2-Sylow subgroups: $\{1,a\}$, 
$\{1,ab\}$, and $\{1,ba\}$, and one 3-Sylow subgroup: $\{1,b,b^2\}$.

Let $p$ be the characteristic of the field $k$; we allow the possibility that $p=0$.

\begin{case}[$\alpha$]
If $p\ne 2,3$, every $\gring k{\Sigma_3}$-module is
projective, by Theorem \eqref{maschke}.
\end{case}

\begin{case}[$\beta$]
Suppose $p=3$.  By \eqref{projectives:sylows}, the rank of each
finite projective $\gring kG$-module is divisible by 3, since the latter is the order of
the 3-Sylow subgroup. Moreover, (\ref{blocks}.d) implies that the number of indecomposable
projectives equals the number of simple modules, and the latter is at least $2$, for
example, by Exercise \eqref{local=pgroup}.  These lead us to the conclusion that there are
exactly two indecomposable projectives, each having rank $3$.

One way to construct them is as follows: Let $H=\{1,a\}$, a 2-Sylow subgroup of
$\Sigma_3$.  There are two nonisomorphic $\gring kH$-module structures on $k$: the
trivial one, given by the augmentation map, and the one defined by character $\sigma\col
H\to k$ with $\sigma(a)=-1$; denote the latter ${}^\sigma\!k$. Plainly, both these $\gring
kH$-modules are simple and hence, by Theorem \eqref{maschke}, projective.  Consequently,
base change along the canonical inclusion $\gring kH\to \gring k{\Sigma_3}$ gives us two
projective $\gring k{\Sigma_3}$-modules,
\[
\gring k{\Sigma_3}\otimes_{\gring kH} k\quad\text{and}\quad
\gring k{\Sigma_3}\otimes_{\gring kH} {}^\sigma\!k.
\]
They both have rank 3. I leave it to you to verify that they are not isomorphic.
Hint: calculate the $\Sigma_3$-invariants.
\end{case}

\begin{case}[$\gamma$]
The situation gets even more interesting when $p=2$.
I claim that there are two indecomposable projective $\gring kG$-modules, of ranks $2$ and
$4$, when $x^2+x+1$ is irreducible in $k$, and three of them, each of rank $2$, otherwise.

Indeed, let $H=\{1,b,b^2\}$; this is a cyclic group of order 3. Hence, by
\eqref{projectives:cyclic}, when $x^2+x+1$ is irreducible in $k[x]$, the group algebra
$\gring kH$ has $2$ (nonisomorphic) simple modules, of ranks $1$ and $2$, and when
$x^2+x+1$ factors in $k[x]$, there are $3$ simple modules, each of rank $1$.  As the
characteristic of $k$ does not divide $|H|$, all these simple modules are projective, so
base change along the inclusion $\gring kH\subset \gring k{\Sigma_3}$ gives rise to the
desired number of projective modules, and of the right ranks, over $\gring kG$.  Note
that, by \eqref{projectives:sylows}, projective modules of rank $2$ are indecomposable.
Thus, to be sure that these are the projectives one seeks, one has to verify that in the
former case the rank $4$ module is indecomposable, and in the latter that the three rank
$2$ modules are nonisomorphic.  Once again, I will let you check this.
\end{case}
\end{chunk}

\section{Cohomology of Supplemented Algebras}
\label{Cohomology of supplemented algebras} 
This section collects basic facts concerning 
the cohomology of supplemented algebras. To
begin with, recall that in the language of Cartan and Eilenberg
\citeyear{Cartan/Eilenberg:1956a} a 
\textit{supplemented $k$-algebra} is a $k$-algebra $R$ with unit $\eta\col k\to R$ and an
augmentation $\eps\col R\to k$ such that $\eps\circ\eta$ is the identity on $k$.

Group algebras are supplemented, but there are many more examples. Take, for instance, any
positively (or negatively) graded $k$-algebra with degree $0$ component equal to $k$.  Or,
for that matter, take the power series ring $k[\![x_1,\cdots,x_n]\!]$, with $\eta$ the
canonical inclusion, and $\eps$ the evaluation at $0$. More generally, thanks to Cohen's
Structure Theorem, if a complete commutative local ring $R$, with residue field $k$,
contains a field, then $R$ is a supplemented $k$-algebra.

Let $R$ be a supplemented $k$-algebra, and view $k$ as an $R$-module via the augmentation.
Let $M$ be a (left) $R$-module. The \textit{cohomology of $R$ with coefficients in $M$} is
the graded $k$-vector space $\Ext*RkM$. The cohomology of $R$ with coefficients in $k$,
that is to say, $\Ext*Rkk$, is usually called the \textit{cohomology of $R$}.{\looseness=-1\par}

The $k$-vector space structure on $\Ext*Rkk$ can be enriched to that of a supplemented
$k$-algebra, and then $\Ext*RkM$ can be made into a \textit{right} module over it.  There
are two ways to introduce these structures: via Yoneda splicing and via compositions.
They yield the same result, up to a sign; see \eqref{product:yoneda}. I have opted for composition
products because it is this description that I use to calculate group cohomology in the
sequel.

\begin{chunk}{\obf Composition products.}
\label{product:composition}
Let $P$ be a projective resolution of $k$.  Composition endows the complex of
$k$-vector spaces $\Hom RPP$ with a product structure, and this product is compatible with
the differential, in the sense that, for every pair of homogenous elements $f,g$ in $\Hom
RPP$, one has
\[
\dd(fg) = \dd(f)g+\sgn ff\dd(g).
\]
In other words, $\Hom RPP$ is a differential graded algebra (DGA). One often refers to this as
the \textit{endomorphism DGA} of $P$. It is not hard to verify that the multiplication of
$\Hom RPP$ descends to homology, that is to say, to $\Ext*Rkk$. This is the
\textit{composition product} on cohomology, and it makes it a graded $k$-algebra. It is even
supplemented, since $\Ext0Rkk=k$.

Let $F$ be a projective resolution of $M$. The endomorphism DGA $\Hom RPP$ acts on the
complex $\Hom RPF$ via composition on the right, and, once again, this action is
compatible with the differentials.  Thus, $\Hom RPF$ becomes a DG \textit{right} module over
$\Hom RPP$. These structures are inherited by the corresponding homology vector spaces;
thus does $\Ext*RkM$ become a right $\Ext*Rkk$-module.

One has to check that the composition products  defined do not depend on the choice of
resolutions; \cite[(7.2)]{Bo} justifies this, and much more.
\end{chunk}

\begin{remark}
\label{product:yoneda}
As mentioned before, one can introduce products on $\Ext*Rkk$ also via Yoneda
multiplication, and, \textit{up to a sign}, this agrees with the composition product; 
\cite[(7.4)]{Bo} has a careful treatment of these issues. The upshot is that one can set
up an isomorphism of $k$-algebras between the Yoneda Ext-algebra and Ext-algebra with
composition products.  Thus, one has the freedom to use either structure, as long as it
is done consistently.
\end{remark}

\begin{chunk}{\obf Graded-commutativity.}
  Let $E$ be a graded algebra. Elements $x$ and $y$ in $E$ are said to commute, in the
  graded sense of the word, if
\[
xy = \sign xy\, yx.
\]
If every pair of its elements commute, $E$ is said to
be graded-commutative.  When $E$ is concentrated in degree $0$ or in even degrees, it is
graded-commutative precisely when it is commutative in the usual sense.

An exterior algebra on a finite-dimensional vector space sitting in odd degrees is another
important example of a graded-commutative algebra. More generally, given a graded
vector space $V$, with $V_i=0$ for $i<0$, the tensor product of the symmetric
algebra on $V_\even$ and exterior algebra on $V_\odd$, that is to say, the $k$-algebra
\[
\sym(V_\even)\otimes_k\bigwedge V_\odd,
\]
is graded-commutative.  If the characteristic of $k$ happens to be $2$, then $\sym(V)$ is
also graded-commutative even when $V_\odd\ne 0$.  This fails in odd characteristics, the
point being that, in a graded-commutative algebra, for an element $x$ of odd degree, $x^2
= -x^2$, so that $x^2=0$ when $2$ is invertible in $E$.

A graded-commutative algebra with the property that $x^2=0$ whenever the degree of
$x$ is odd is said to be \textit{strictly graded-commutative}.  An exterior algebra (with
generators in odd degrees) is one example. Here is one more, closer to home: for
a homomorphism of commutative rings $R\to S$, the graded $S$-module $\Tor *RSS$ is strictly
graded-commutative, with the pitchfork product (homology product) defined by Cartan and
Eilenberg; see \cite[VIII\,\S2]{Mc:hom}.
 \end{chunk}

\begin{chunk}{\obf Functoriality.}
\label{functoriality}
  The product in cohomology is functorial, in that, given a homomorphism of supplemented
  $k$-algebras $\vf\col R\to R'$, the induced map of graded $k$-vector spaces
\[
\Ext*{\vf}kk \col \Ext{*}{R'}kk \to \Ext*Rkk
\]
is a homomorphism of supplemented $k$-algebras.

Now let $R$ and $S$ be supplemented $k$-algebras. The tensor product $R\otimes_kS$ is also
a supplemented $k$-algebra, and the canonical maps
\[
R\alto{1\otimes\eps^S}R\otimes_kS \arto{\eps^R\otimes1}S
 \]
 respect this structure. By functoriality of products, the diagram above induces
 homomorphisms of supplemented $k$-algebras
\[
\Ext*Rkk \arto{\Ext*{1\otimes\eps^S}kk} \Ext*{R\otimes_kS}kk
\alto{\Ext*{\eps^R\otimes1}kk}\Ext*Skk .
\]
It is not hard to check that the images of these maps commute, in the graded sense, so one
has a diagram of supplemented $k$-algebras:
\begin{gather*}
  \xymatrixrowsep{2pc} \xymatrixcolsep{3pc} \xymatrix{ \Ext*Rkk\ar[r]^{\Ext
      *{\id\otimes\eps^S}kk}&\Ext *{R\otimes_kS}kk
    &\Ext*Skk\ar[l]_{\Ext *{\eps^R\otimes\id}kk} \\
    &\Ext*Rkk\otimes_k\Ext *Skk\ar[u]\ar@{<-}[ul]^{\id\otimes
      1}\ar@{<-}[ur]_{1\otimes\id}}\tag{$*$}
\end{gather*}
I should point out that the tensor product on the lower row is the \textit{graded} tensor
product and the multiplication on it is defined accordingly, that is,
\[
(r\otimes_ks)\cdot(r'\otimes_ks') = \sign s{r'} (rr'\otimes_kss').
\]
Under suitable finiteness hypotheses\emdash for example, if $R$ and $S$ are noetherian\emdash the
vertical map in ($*$) is bijective.  However, this is not of
importance to us.
\end{chunk}


\subsection*{The cohomology of Hopf algebras.}
The remainder of this section deals with the cohomology of Hopf algebras. So let $H$ be a
Hopf algebra, with diagonal $\Delta$ and augmentation $\eps$; see \eqref{hopf algebras}.
The main example to keep in mind is the case when $H$ is the group algebra of a group,
with the diagonal defined in \eqref{gring:diagonal}.

One crucial property of the cohomology algebra of $H$, which distinguishes it from the
cohomology of an arbitrary supplemented algebra, is the following.

\begin{proposition}\label{product:commutative}
  The cohomology algebra $\Ext *Hkk$ is graded-commutative.
 \end{proposition}
 
 Note that $H$ is not assumed to be cocommutative.  This is a striking result,
 and its proof is based on the diagram of $k$-algebra homomorphisms
\begin{equation}
\label{cups}
\Ext*Hkk\otimes_k\Ext*Hkk\to \Ext*{H\otimes_kH}kk\arto{\Ext*{\Delta}kk} \Ext *Hkk,
  \end{equation}
  where the one on the left is the vertical map in (\ref{functoriality}.1), with $R$ and $S$
  equal to $H$, and the one on the right is induced by the diagonal homomorphism.

\begin{proposition}\label{product:delta}
  The composition of homomorphisms in \textit{(\ref{product:commutative}.\ref{cups})} is the
  product map; that is to say, $(x\otimes_ky)\mapsto xy$ for $x$ and $y$ in
  \textit{$\Ext*Hkk$}.
  
  In particular, the product map of \textit{$\Ext*Hkk$} is a homomorphism of $k$-algebras.
\regulardigits\hfuzz.2pt\par
\end{proposition}
\begin{proof}
  The diagram in question expands to the following commutative diagram of homomorphisms of
  $k$-algebras, where the lower half is obtained from (\ref{functoriality}.1), the upper
  half is induced by property (c) of Hopf algebras\emdash see
  \eqref{hopf algebras}\emdash to the effect
  that $\eps$ is a co-unit for the diagonal.
\[
\xymatrixrowsep{2pc} \xymatrixcolsep{3pc} \xymatrix{
  &\Ext *Hkk\ar@{<-}[d]^{\Ext*{\Delta}kk}\ar@{<-}[dl]_{\id}\ar@{<-}[dr]^{\id}& \\
  \Ext*Hkk\ar[r]_{\Ext *{\id\otimes\eps}kk}&\Ext *{H\otimes_kH}kk
  &\Ext*Hkk\ar[l]^{\Ext *{\eps\otimes\id}kk} \\
  &\Ext*Hkk\otimes_k\Ext *Hkk\ar[u]\ar@{<-}[ul]^{\id\otimes 1}\ar@{<-}[ur]_{1\otimes\id}}
\]
Let $x$ and $y$ be elements in $\Ext*Hkk$.  The element $x$ goes to $x\otimes_k1$ under
the map heading southeast, and to $x$ under the map heading northeast. The
commutativity of the diagram thus implies that $x\otimes_k1\mapsto x$ under the composed
vertical map. A similar diagram chase reveals that $1\otimes_ky\mapsto y$. Since the
vertical maps are homomorphisms of $k$-algebras, one has 
\[
x\otimes_ky = (x\otimes_k1)\cdot(1\otimes_ky) \mapsto xy.
\]
This is the conclusion we seek.
\end{proof}

The proof of Proposition \eqref{product:commutative} uses also the following elementary
exercise, of which there are versions for groups, for coalgebras, etc.

\begin{exercise}
\label{above}
A graded $k$-algebra $R$ is graded-commutative precisely when the product map
  $R\otimes_k R\to R$ with $r\otimes s\mapsto rs$ is a homomorphism of rings.
\end{exercise}

Now one can prove that the cohomology algebra is graded-commutative.

\begin{proof}[Proof of Proposition \eqref{product:commutative}]
  By the preceding proposition, the product map
$\Ext*Hkk \otimes_k \Ext*Hkk\to \Ext *Hkk$ given by
$x\otimes_ky\mapsto xy$ is a homomorphism of rings (for a general
algebra it is only $k$-linear). 
To complete the proof one has to do Exercise~\eqref{above}.
 \end{proof}

\section{Group Cohomology}
In this section we return to group algebras.

\begin{chunk}{\obf Cohomology.}
  Let $G$ be a group and let $M$ be a $\gring kG$-module. Recall that $\gring kG$
  is a supplemented algebra. The \textit{cohomology of $\,G$ with coefficients in $M$} is
  the graded $k$-vector space
\[
\gc *GM = \Ext *{\gring kG}kM.
\]
There is no ambiguity concerning the field $k$ since $\Ext *{\gring kG}kM$ is
isomorphic to $\Ext *{\gring {\BZ}G}{\BZ}M$; see \cite[(1.1)]{Evens:1961a}. The \textit{cohomology of\,
  $G$} is $\gc *Gk$.

Standard properties of Ext-modules carry over to the situation on hand.  For instance,
each short exact sequence of $\gring kG$-modules $0\to L\to M\to N\to 0$ engenders a long
exact sequence of $k$-vector spaces
\[
0\to \gc 0GL \to \gc 0GM \to \gc 0GN \to \gc 1GL \to \gc 1GM \to \cdots .
\]

Note that $\gc nG-=0$ for $n\geq 1$ if and only if $k$ is projective.  Therefore, one has
the following cohomological avatar of Maschke's theorem \eqref{maschke}:

\begin{theorem}
\label{gc:maschke}
Let $G$ be a finite group.  Then $\gc nG-=0$ for each integer $n\geq 1$ if and only if the
characteristic of $k$ is coprime to $|G|$. \hfill$\square$
\end{theorem}

As is typical in homological algebra, low degree cohomology modules have nice
interpretations.  For a start, $\Ext 0{\gring kG}kM=\Hom{\gring kG}kM$, so
\eqref{homs} yields
\[
\gc 0GM = M^G .
\]
Thus, one can view the functors $\gc nG-$ as the derived functors of invariants.

The degree $1$ component of $\gc *GM$ is also pretty down to earth.  Recall that a map
$\theta\col G\to M$ is said to be a \textit{derivation}, or a \textit{crossed homomorphism},
if it satisfies the Leibniz formula: $\theta(gh)=\theta(g)+ g\theta(h)$, for every $g,h$
in $G$.  The asymmetry in the Leibniz rule is explained when one views $M$, which is
\textit{a priori} only a left $\gring kG$-module, as a $\gring kG$-bimodule with trivial
right action: $m\cdot g = m$. Using the $k$-vector space structure on $M$ one can add
derivations, and multiply them with elements in $k$, so they form a $k$-vector space; this
is denoted $\der GM$. This vector space interests us because of the following
\begin{lemma}
\label{derivations}  
The $k$-vector spaces $\Hom{\gring kG}{\aug G}M$ and $\der GM$ are isomorphic via the maps
\[
\begin{gathered}
{\begin{aligned}
  \Hom{\gring kG}{\aug G}M &\to \der GM  \\               
  \qquad\alpha &\mapsto \big(g\mapsto \alpha(g-1)\big),\\
 \end{aligned}}\qquad
{\begin{aligned}
\der GM & \to\Hom{\gring kG}{\aug G}M \\
\qquad \theta &\mapsto \big(g-1\mapsto \theta(g)\big).
\end{aligned}}
\end{gathered}
\]
\end{lemma}
The proof is an elegant computation and is best rediscovered on one's own.  As to its
bearing on $\gc 1GM$: applying $\Hom{\gring kG}-M$ to the exact sequence
\[
0\to \aug G \to \gring kG\to k\to 0
\]
of $\gring kG$-modules leads to the exact sequence of $k$-vector spaces
\[
0\to M^G \to M \to \der GM\to \gc 1GM\to 0.
\]
In this sequence, each $m\in M$ maps to a derivation: $g\mapsto (g{-}1)m$; these are the
\textit{inner derivations} from $G$ to $M$, and their set is denoted
by $\inn GM$.  Thus,
\[
\gc 1GM = \der GM/\inn GM.
\]
Let us specialize to the case when $M=k$. The Leibniz rule for a derivation
$\theta\col G\to k$ then reads: $\theta(gh) = \theta(g) + \theta(h)$, so $\der Gk$ coincides
with group homomorphisms from $G$ to $k$. Moreover, every inner derivation from
$G$ to $k$ is trivial.  The long and short of this discussion is that $\gc 1Gk$ is
precisely the set of additive characters from $G$ to $k$.

There are other descriptions, some of a more group theoretic
flavour, for $\gc 1GM$; for those the reader may look in \cite{Benson:1991a}.
 \end{chunk}
 
 The discussion in Section \ref{Cohomology of supplemented algebras} on products on
 cohomology applies in the special case of the cohomology of group algebras. In
 particular, since $\gring kG$ is a Hopf algebra, Proposition \eqref{product:commutative}
 specializes thus:
\begin{theorem}
 The cohomology algebra $\gc *Gk$ is graded-commutative.\hfill$\square$
\end{theorem}

\begin{chunk}{\obf K\"unneth formula.}
\label{kunneth}
Let $G_1$ and $G_2$ be groups. Specializing (\ref{functoriality}.1) to the case where
$R=\gring k{G_1}$ and $S=\gring k{G_2}$, one obtains a homomorphism of $k$-algebras
\[
\gc *{G_1}k\otimes_k\gc *{G_2}k \to \gc *{G_1\times G_2}k.
\]
This map is bijective whenever the group algebras are noetherian. This is
the case when, for example, $G_i$ is finite, or finitely generated and
abelian.
\end{chunk}

\begin{chunk}{\obf Resolutions.}
\label{resolutions}
If one wants to compute cohomology from first principles, one has to first obtain a
projective resolution of $k$ over $\gring kG$. In this regard, it is of interest to get as
economical a resolution as possible. Fortunately, any finitely generated module over $\gring
kG$ has a minimal projective resolution; we discussed this point already in
\eqref{covers}; unfortunately, writing down this minimal resolution is a challenge. In
this the situation over group algebras is similar to that over commutative local rings.
What is more difficult is calculating products from these minimal resolutions.

There is  a canonical resolution for $k$ over $\gring kG$ called the \textit{Bar
  resolution}; while it is never minimal, it has the merit that there is a simple formula
for calculating the product of cohomology classes. The are many readable sources for
this, such as \cite[(3.4)]{Benson:1991a}, \cite[(2.3)]{Evens:1991a},
  and \cite[IV\,\S5]{Mc:hom}, so I will not 
reproduce the details here.
\end{chunk}

\section{Finite Generation of the Cohomology Algebra}

In the preceding section, we noted that the cohomology algebra of a finite group is
graded-commutative. From this, the natural progression is to the following theorem,
contained in \cite{Evens:1991a}, \cite{Golod:1959a}, and
\cite{Venkov:1959a}.

\begin{theorem}
\label{fingen}
Let $G$ be a finite group. The $k$-algebra $\gc *Gk$ is finitely generated, and hence
noetherian. \hfill$\square$
\end{theorem}  

This result, and its analogues for other types of groups, is the starting point of
Benson's article \citeyear{Be:msri}; see the discussion in Section 4
of it. There are many ways of 
proving Theorem~\eqref{fingen}, some more topological than others; one
that is entirely algebraic is given in
\cite[(7.4)]{Evens:1961a}.

In this section I prove the theorem in some special cases.  But first:

\begin{ramble}
  Theorem \eqref{fingen} has an analogue in commutative algebra: Gulliksen \citeyear{Gu}
  proves that when a commutative local ring $R$, with residue field
  $k$, is a complete
  intersection, the cohomology algebra $\Ext{*}Rkk$ is noetherian. There is
  a perfect converse: B{\o}gvad and Halperin \citeyear{Bogvad/Halperin:1986a} have proved that if the $k$-algebra
  $\Ext{*}Rkk$ is noetherian, then $R$ must be complete intersection.
  
  There are deep connections between the cohomology of modules over complete intersections
  and over group algebras.  This is best illustrated by the theory of support varieties.
  In group cohomology it was initiated by Quillen \citeyear{Quillen:1971b,Quillen:1971c}, and developed in depth by
  Benson and Carlson, among others; see \cite{Benson:1991b} for a
  systematic introduction.  In commutative algebra, support varieties
  were introduced by Avramov \citeyear{Av:vcd}; see also \cite{AvBu}.
  
  As always, there are important distinctions between the two contexts. For example, the
  cohomology algebra of a complete intersection ring is generated by its elements of
  degree $1$ and $2$, which need not be the case with group algebras. More importantly,
  once the defining relations of the complete intersection are given, one can write down
  the cohomology algebra; the prescription for doing so is given in \cite{So}.
  Computing group cohomology is an entirely different cup of tea.  Look up \cite{Carlson:2001a} for
  more information on the computational aspects of this topic.
\end{ramble}

Now I describe the cohomology algebra of finitely generated abelian groups.  In
this case, the group algebra is a complete intersection\emdash see \eqref{gring:abelian}\emdash so
one may view the results below as being about commutative rings or about finite groups.

\begin{proposition}
\label{gc:zn}
For each positive integer $n$, the cohomology of $\BZ^n$ is the exterior algebra on an
$n$-dimensional vector space concentrated in degree $1$.
\end{proposition}
\begin{proof}
  As noted in \eqref{gring:cyclic}, the group algebra of $\BZ$ is $k[x^{\pm1}]$, with
  augmentation defined by $\eps(x)=1$.  The augmentation ideal is generated by $x-1$, and
  since this element is regular, the Koszul complex
\[
0\to k[x^{\pm1}] \arto{x-1}k[x^{\pm1}] \to 0,
\]
is a free resolution of $k$. Applying $\Hom {k[x^{\pm1}]}-k$ yields the complex with
trivial differentials: $0\to k\to k\to 0$, and situated in cohomological degrees $0$ and
$1$. Thus, $\gc 0{\BZ}k=k=\gc 1{\BZ}k$.  Moreover, $\gc 1{\BZ}k\cdot \gc 1{\BZ}k=0$, by
degree considerations, so that the cohomology algebra is the exterior algebra $\wedge_k\,
k$, where the generator for $k$ sits in degree $1$.

For $\BZ^n$, one uses the K\"unneth formula \eqref{kunneth} to calculate group cohomology:
\[
\gc *{\BZ^n}k = \gc *{\BZ}k^{\otimes n} = \wedge_k\, k^n,
\]
where the generators of $k^n$ are all in (cohomological) degree $1$.
\end{proof}

The next proposition computes the cohomology of cyclic $p$-groups.  It turns out that one
gets the same answer for all but one of them; the odd man out is the group of order two.

\begin{proposition}
\label{gc:pgroups}
Let $k$ be a field of characteristic $p$, and let $G=\BZ/p^e\BZ$, for some integer
$e\geq1$.
\begin{enumerate}
\item[{\rm(i)}] When $p=2$ and $e=1$, $\gc *Gk=\sym(k\dual e_1)$, with $|\dual e_1|=1$.
\item[{\rm(ii)}] Otherwise $\gc *Gk = \bigwedge  (k\dual e_1)\otimes_k\sym(k\dual e_2)$, with $|\dual
  e_1|=1$ and $|\dual e_2|=2$.
 \end{enumerate}
 \end{proposition}
\begin{proof} 
  The group algebra of $G$ is $k[x]/(x^{p^e}-1)$, and its augmentation
  ideal is $(x-1)$. Note that $x^{p^e}-1 = (x-1)^{p^e}$, so the substitution $y=x-1$
  presents the group algebra in the more psychologically comforting, to this commutative
  algebraist, form $k[y]/(y^{p^e})$. Write $R$ for this algebra; it is a $0$-dimensional
  hypersurface ring\emdash the simplest example of a complete intersection\emdash with socle
  generated by the element $y^{p^e-1}$. The $R$-module $k$ has minimal free resolution
\[
P\col \cdots \to Re_3 \arto y Re_2\arto{y^{p^e-1}} Re_1\arto y Re_0\to 0.
\] 
This is an elementary instance of the periodic minimal free resolution, of period 2, of
the residue field of hypersurfaces constructed by Tate \citeyear{Ta}; see also 
\cite{Ei}. Applying $\Hom R-k$ to the resolution above results in the complex
\[
\Hom RPk \col 0\to k\dual e_0\arto{0}k{\dual e_1}\arto{0}k{\dual e_2}\arto{0}k{\dual
  e_3}\arto{0} \cdots
\]
Thus, one obtains  $\gc nGk=k$ for each integer $n\geq 0$.

\subsubsection*{Multiplicative structure.}
Next we calculate the products in group cohomology, and for this I propose to use
compositions in $\Hom RPP$; see \eqref{product:composition}. More precisely: since $P$ is
a complex of free modules, the canonical map
\[
\Hom RP\eps\col \Hom RPP\to \Hom RPk
\]
is an isomorphism in homology. Given two cycles in $\Hom RPk$, I will lift them to cycles
in $\Hom RPP$, compose them there, and then push down the resultant cycle to $\Hom RPk$;
this is their product.

For example, the cycle $\dual e_1$ of degree $-1$ lifts to the cycle $\alpha$ in $\Hom RPP$
given by
\[
\xymatrixrowsep{2pc} \xymatrixcolsep{2.5pc} \xymatrix{ \cdots \ar[r]\ar[rd]^1 &
  Re_4\ar[r]^{y^{p^e-1}}\ar[rd]^{-y^{p^e-2}}& Re_3\ar[r]^y\ar[rd]^1&
  Re_2\ar[r]^{y^{p^e-1}}\ar[rd]^{-y^{p^e-2}} &
  Re_1\ar[r]^y\ar[rd]^1 &  Re_0\ar[r] & 0 \\
  \cdots \ar[r]& Re_4\ar[r]_{y^{p^e-1}}& Re_3\ar[r]_y & Re_2\ar[r]_{y^{p^e-1}}&
  Re_1\ar[r]_y&Re_0\ar[r] & 0 }
\]
It is a lifting of $\dual e_1$ since $\eps(\alpha(e_1))=1$, and a cycle since
$\dd\alpha=-\alpha\dd$. Similarly, the cycle $\dual e_2$ lifts to the cycle $\beta$ given
by
\[
\xymatrixrowsep{2pc} \xymatrixcolsep{2.5pc} \xymatrix{
  \cdots\ar[r]\ar[rrd]^1&Re_4\ar[r]^{y^{p^e-1}}\ar[rrd]^1&
  Re_3\ar[r]^y\ar[rrd]^1&Re_2\ar[r]^{y^{p^e-1}}\ar[rrd]^1 &
  Re_1\ar[r]^y& Re_0\ar[r] & 0 \\
  \cdots\ar[r] &Re_4\ar[r]_{y^{p^e-1}}&Re_3\ar[r]_y & Re_2\ar[r]_{y^{p^e-1}}&
  Re_1\ar[r]_y&Re_0\ar[r] & 0 }
\]

This is all one needs in order to compute the entire cohomology rings of $G$. As indicated before,
there are two cases to consider.

When $p=2$ and $e=1$, one has $y^{p^e-2}=1$, so that $\eps(\alpha^n(e_n))=1$ for each
positive integer $n$. Therefore, $(\dual e_1)^n=\dual e_n$, and since the $\dual e_n$ form
a basis for the graded $k$-vector space $\gc *Gk$, one obtains $\gc *Gk=k[\dual e_1]$,
as desired.

Suppose that either $p\geq 3$ or $e\geq 2$. In this case
\[
\eps(\alpha^{n+1}(e_{n+1})) = 0, \quad \eps(\beta^n(e_{2n}))=1,\quad
\hbox{and}\quad\eps(\alpha\beta^{n-1}(e_{2n-1}))=1,
\]
for each positive integer $n$.  Passing to $\Hom RPk$, these relations translate to
\[
(\dual e_1)^{n+1} = 0, \quad (\dual e_2)^n=\dual e_{2n}, \quad \dual e_1(\dual
e_2)^{n-1}=\dual e_{2n-1}.
\]
In particular, the homomorphism of $k$-algebras $k[\dual e_1,\dual e_2]\to \gc *Gk$ is
surjective; here, $k[\dual e_1,\dual e_2]$ is the graded-polynomial algebra on $\dual e_1$
and $\dual e_2$, that is to say, it is the tensor product of the exterior algebra on
$\dual e_1$ and the usual polynomial algebra on $\dual e_2$.  This map is also injective:
just compare Hilbert series.

This completes our calculation of the cohomology of cyclic $p$-groups.
\end{proof}

\begin{chunk}{\obf Finitely generated abelian groups.}
\label{cohomology:abelian}
Let the characteristic of $k$ be $p$, and let the group $G$ be finitely generated and abelian.
By the fundamental theorem of finitely generated abelian groups, there are integers $n$
and $e_1,\dots,e_m$, such that
\[ 
G \cong \BZ^n \oplus\frac\BZ{(p^{e_1}\BZ)}\oplus \cdots \oplus\frac
\BZ{(p^{e_m}\BZ)}\oplus G'.
\]
where $G'$ is a finite abelian group whose order is coprime to $p$.  By the K\"unneth
formula \eqref{kunneth}, the group cohomology of $G$ is the $k$-algebra
\[
\gc *{\BZ^n}k\otimes_k \gc *{\BZ/p^{e_1}\BZ}k \otimes_k \cdots \otimes_k \gc
*{\BZ/p^{e_m}\BZ}k\otimes_k\gc *{G'}k.
\]
Note that $\gc *{G'}k=k$, by Theorem \eqref{gc:maschke}; the remaining terms of
the tensor product above are computed by propositions \eqref{gc:zn} and
\eqref{gc:pgroups}. 
\end{chunk}

To give a flavour of the issues that may arise in the nonabelian case, I will 
calculate the cohomology of $\Sigma_3$. This gives me also an excuse to introduce an
important tool in this subject:

\begin{chunk}{\obf The Lyndon--Hochschild--Serre spectral sequence.}
\label{hsss}
Let $G$ be a finite group and $M$ a $\gring kG$-module. Let $N$ be a normal subgroup in
$G$.
  
Via the canonical inclusion of $k$-algebras $\gring kN\subseteq \gring kG$, one can view
$M$ also as an $\gring kN$-module.  Since $N$ is a normal subgroup, the $k$-subspace $M^N$
of $N$-invariant elements of $M$ is stable under multiplication by elements in $G$
(check!) and hence it is a $\gring kG$-submodule of $M$. Furthermore, $\aug N\cdot M^N=0$,
so that $M^N$ has the structure of a module over $\gring kG/\aug N\gring kG$, that is to
say, of a $\gring k{G/N}$-module; see \eqref{gring:functor}. It is clear from the
definitions that $(M^N)^{G/N}=M^G$. In other words, one has an isomorphism of functors
\[
\Hom {\gring k{G/N}}k{\Hom{\gring kN}k-}\cong \Hom{\gring kG}k- .
\]
The functor on the left is the composition of two functors: $\Hom{\gring kN}k-$ and
$\Hom{\gring k{G/N}}k-$. Thus standard homological algebra provides us with a spectral
sequence that converges to its composition, that is to say, to $\gc{*}GM$.  In our case,
the spectral sequence sits in the first quadrant and has second page
\[
\EC 2pq = \gc p{G/N}{\gc qNM} 
\]
and
differential
\[ \dd^{p,q}_r\col \EC rpq \to
\EC r{p+r}{q-r+1}.
\]
This is the \textit{Lyndon--Hochschild--Serre spectral sequence} associated to $N$.
\end{chunk}

Here are two scenarios where the spectral sequence collapses.

\begin{chunk}
\label{hsss:special1}
Suppose the characteristic of $k$ does not divide $[G:N]$, the index of $N$ in $G$.
In this case, $\gc p{G/N}-=0$ for $p\geq 1$, by Maschke's theorem \eqref{gc:maschke}, so
that the spectral sequence in \eqref{hsss} collapses to yield an isomorphism
\[
\gc *GM \cong \gc 0{G/N}{\gc *NM} = \gc *NM^{G/N}.
\]
In particular, with $M=k$, one obtains that $\gc *Gk\cong \gc *Nk^{G/N}$; this isomorphism
is compatible with the multiplicative structures. Note that the object on the right is the
ring of invariants of the action of $G/N$ on the group cohomology of $N$.  Thus does
invariant theory resurface in group cohomology.
\end{chunk}

\begin{chunk}
\label{hsss:special2}
Suppose the characteristic of $k$ does not divide $|N|$. Then $\gc qNM=0$ for $q\geq
1$, and once again the spectral sequence collapses to yield an isomorphism
\[
\gc *GM \cong \gc *{G/N}{M^N}.
\]
The special case $M=k$ reads $\gc *Gk = \gc *{G/N}k$.
\end{chunk}

As an application we calculate the cohomology of $\Sigma_3$:

\begin{chunk}{\obf The symmetric group on three elements.}
  In the notation in \eqref{projectives:example}, set $N=\{1,b,b^2\}$; this is a normal
  subgroup of $\Sigma_3$, and the quotient group $\Sigma_3/N$ is (isomorphic to)
  $\BZ/2\BZ$. We use the Hochschild--Serre spectral sequence generated by $N$ in
  order to calculate the cohomology of $\Sigma_3$.  There are three cases.

\begin{case}[$\alpha$]
  When $p\ne 2,3$, Maschke's theorem \eqref{gc:maschke} yields
\[
\gc n{\Sigma_3}k \cong
       \begin{cases} k & \text{if $n=0,$}\cr 0
         & \text{otherwise.}\end{cases}
\]
\end{case}

\begin{case}[$\gamma$] If $p=2$, then
\[
\gc *{\Sigma_3}k = k[\dual e_1], \where |\dual e_1|=1;
\]
the polynomial ring on the variable $e_1$ of degree $1$.  Indeed, the order of $N$ is $3$,
so \eqref{hsss:special2} yields that $\gc*{\Sigma_3}k=\gc *{\BZ/2\BZ}k$. Proposition
\eqref{gc:pgroups} does the rest.
 \end{case}

\begin{case}[$\beta$] Suppose that $p=3$. One obtains from \eqref{hsss:special1} that
\[
\gc *{\Sigma_3}k = \gc *Nk^{\BZ/2\BZ}.
\]
The group $N$ is cyclic of order $3$, so its cohomology is $k[\dual e_1,\dual e_2]$, with
$|\dual e_1|=1$ and $|\dual e_2|=2$; see Proposition \eqref{gc:pgroups}. The next step is to
compute the ring of invariants.  The action of $y$, the generator of $\BZ/2\BZ$, on $\gc
*Nk$ is compatible with products, so it is determined entirely by its actions on $\dual e_1$
and on $\dual e_2$. I claim that
\[
y(\dual e_1) = -\dual e_1\qnd y(\dual e_2) = -\dual e_2.
\]
Using the description of $\gc 1Nk$ given in \eqref{derivations}, it is easy to verify the
assertion on the left; the one of the right is a little harder. Perhaps the best way to
get this is to observe that the action of $y$ on $\gc *Nk$ is compatible with the
\textit{Bockstein} operator on cohomology and that this takes $e_1^*$ to $e_2^*$; see
\cite[(3.3)]{Evens:1961a}.  At any rate, given this, it is not hard to see that
\[
\gc *{\Sigma_3}k = \bigwedge(k\dual e_1\dual e_2)\otimes_k\sym(k(\dual e_2)^2),
\]
the tensor product of an exterior algebra on an element of degree $3$ and a symmetric
algebra on an element of degree $4$.
  \end{case}
\end{chunk}

\subsection*{Hopf algebras.}

In this article I have indicated at various points that much of the
module theory over group algebras extends to Hopf algebras.  I 
wrap up by mentioning a perfect generalization of Theorem
\eqref{fingen}, due to E.~Friedlander and Suslin
\citeyear{Friedlander/Suslin:1997a}: If a finite-dimensional Hopf
algebra $H$ is cocommutative, its cohomology algebra $\Ext *Hkk$
is finitely generated.

\def\citeauthoryear#1#2#3{#1 #3}
\bibliographystyle{msribib}
\bibliography{book51}
\end{document}